\documentclass{amsart}
\usepackage{amsmath,amssymb}
\usepackage[mathscr]{eucal}

\usepackage[all]{xy}

\theoremstyle{plain}
\newtheorem{theorem}{Theorem}[section]

\newtheorem{lemma}[theorem]{Lemma}
\newtheorem{corollary}[theorem]{Corollary}

\theoremstyle{definition}
\newtheorem{definition}[theorem]{Definition}
\newtheorem{example}[theorem]{Example}

\theoremstyle{remark}
\newtheorem{remark}[theorem]{Remark}

\newcommand{\secref}[1]{Section~\ref{#1}}
\newcommand{\thmref}[1]{Theorem~\ref{#1}}

\newcommand{\lemref}[1]{Lemma~\ref{#1}}
\newcommand{\corref}[1]{Corollary~\ref{#1}}

\newcommand{\defref}[1]{Definition~\ref{#1}}

\def\Longrightarrow{- \! -  \! \! \! \longrightarrow}

\def\Z{{\mathbb Z}}
\def\Q{{\mathbb Q}}
\def\C{{\mathbb C}}
\def\L{{\mathbb L}}
\def\map{\mathrm{map}}

\def\ad{\mathrm{ad}}
\def\Der{\mathrm{Der}}

\def\Rel{\mathrm{Rel}}
\def\ker{\mathrm{ker}}
\def\im{\mathrm{im}}

\def\QL{\mbox{\boldmath$\mathscr{L}$} }

\begin{document}

\title[Quillen Models and Adjoint Maps]{Rationalized Evaluation Subgroups of a Map II: Quillen Models and Adjoint Maps}

\author{Gregory  Lupton}

\address{Department of Mathematics,
          Cleveland State University,
          Cleveland OH 44115}

\email{G.Lupton@csuohio.edu}

\author{Samuel Bruce Smith}

\address{Department of Mathematics,
  Saint Joseph's University,
  Philadelphia, PA 19131}

\email{smith@sju.edu}

\date{\today}



\begin{abstract}
Let $\omega\colon \map(X, Y;f) \to Y$ denote
a general evaluation
fibration.  Working in the setting of
rational homotopy theory via
differential graded Lie algebras, we
identify the long exact
sequence induced on rational homotopy groups
by $\omega$ in terms of
(generalized) derivation spaces and adjoint
maps. As a
consequence, we obtain a unified description
of the rational
homotopy theory of function spaces, at the
level of rational
homotopy groups, in terms of derivations of
Quillen models and
adjoints.  In particular, as a natural
extension of a result of
Tanr\'{e}, we identify the rationalization
of the evaluation
subgroups of a map $f\colon X \to Y$ in this
setting.  As
applications, we consider a generalization
of a question of
Gottlieb, within the context of rational
homotopy theory. We also
identify the rationalization of the $G$-sequence of $f$ and make
explicit computations of the homology of
this sequence.  In a
separate result of independent interest, we
give an explicit
Quillen minimal model of a product $A\times X$, in the case in
which $A$ is a rational co-$H$-space.
\end{abstract}

\maketitle

\section{Introduction}\label{section:intro}

Let $f\colon X \rightarrow Y$ be a based map
of connected spaces.
Let $\map(X, Y)$ denote the space of unbased
maps from $X$ to $Y$
and  $\map(X, Y;f)$  the path component
containing $f$. Evaluation
at the basepoint of $X$ determines a
fibration $\omega \colon
\map(X, Y;f) \rightarrow Y$. The $n$th
\emph{evaluation subgroup of $f$} is then defined to be the subgroup
$$ G_n(Y, X;f) = \hbox{Image} \{ \omega_\#
\colon
\pi_n(\map(X, Y;f)) \to \pi_n(Y) \}$$
of $\pi_n(Y)$. The $n$th Gottlieb group
$G_n(X)$ of a space $X$
occurs as the special case $X=Y$ and $f=1$.
Because $\omega\colon
\map(X,X;1) \rightarrow X$ can be identified
with the connecting
map of the universal fibration for
fibrations with fibre $X$, the
Gottlieb group is an important universal
object for fibrations
with fibre $X$. The structure of the
Gottlieb group and its role
in homotopy theory represents a broad area
of research having its
beginnings, of course, with the papers of
Gottlieb \cite{Got1,
Got2,  Go1} (see \cite{Oprea} for a recent
survey and further
references). A fundamental obstacle to
computing Gottlieb groups
is their lack of functorality, in that
generally
$f_\#\big(G_*(X)\big) \not\subseteq G_*(Y)$
for a map $f\colon X
\to Y$.  By widening our perspective so as
to include the
evaluation subgroups of a map, we remedy
this situation somewhat,
in that $f\colon X \to Y$ always induces a
map $f_\#\colon G_n(X)
\to G_n(Y, X;f)$.

In a previous paper \cite{L-S}, we have
developed a basic
framework within which rational homotopy
groups of function spaces
and related topics---including rationalized
evaluation
subgroups--- may be studied.  All results in
that paper were
developed from the Sullivan minimal model
point of view, that is,
using DG algebras. The current paper is
intended both as a
complement to and a continuation of the
earlier work. It is
complementary, in that all results here, and
the basic framework
that we establish, are developed from the
Quillen minimal model
point of view, that is, using DG Lie
algebras.  It continues the
earlier work in that we present some
different developments of the
basic results. We refer to \cite{L-S} for a
general introduction
to the themes of this paper, and for a
survey of the exisiting
results in the area.

In \cite[Th.2.1]{L-S}, we identified the
rationalization of the
induced homomorphism $\omega_\#\colon
\pi_n(\map(X, Y;f)) \to
\pi_n(Y)$ in terms of the homology of
derivation spaces of the
Sullivan models of the spaces $X$ and $Y$.
Our result extends, to
a connected component of a general function
space, the isomorphism
$$\pi_n(\map(X,X;1)) \otimes \Q \cong
H_n(\Der_*(\mathcal{M}_X)),$$
which was first observed by Sullivan in
\cite{Su}. Here
$\Der_*(\mathcal{M}_X)$ is the differential
graded vector space of
degree-lowering derivations of the Sullivan
model of $X$. Our
results in \cite{L-S} provide a general
framework for studying
certain long exact sequences of homotopy
groups for function space
components within the context of Sullivan
minimal models. In
particular, they allow for an extension of
the F\'{e}lix-Halperin
description of the rational Gottlieb groups
(see \cite[Th.3]{F-H})
to the more general setting of rational
evaluation subgroups of a
map.

The developments of this paper proceed in an
analogous way to
those of \cite{L-S}.  As a special case of
\thmref{thm:isos Phi
and Phi*}, we obtain an isomorphism
$$\pi_n(\map_*(X,X;1)) \otimes \Q \cong
H_n(\Der_*(\mathcal{L}_X)),$$
where $(\Der_*(\mathcal{L}_X))$ is the
differential graded vector
space of degree-raising derivations of the
Quillen model of $X$.
More generally, the results of
\secref{section:main result} lead
to a complete picture, within the framework
of Quillen models, of
various long exact sequences of rational
homotopy groups of
function space components.   In particular,
in
\thmref{theorem:evaluation group} we  extend
Tanr\'{e}'s
description of the rationalized Gottlieb
group
\cite[VII.4.(10)]{Tan83} to a description of
the rationalized
evaluation subgroup of a map.

We pursue two applications of these results
here. We first
consider a generalization  of a question of
Gottlieb, in the
context of rational homotopy theory. It is
well-known that
Gottlieb elements in $G_*(X)$ have vanishing
Whitehead product
with all elements  of $\pi_*(X)$. Let $[ \ \
, \ \ ]_w$ denote the
Whitehead bracket in $\pi_*(X)$ and let
$P_*(X)$ denote the
subgroup of $\pi_*(X)$ consisting of
homotopy elements with
vanishing Whitehead product with all
elements of $\pi_*(X)$---the
so-called \emph{Whitehead center} of
$\pi_*(X)$. Gottlieb's
question asks about the difference between
$P_*(X)$ and its
subgroup $G_*(X)$.  More generally, given a
map $f\colon X \to Y$
set
$$P_n(Y,X;f) = \{ \alpha \in \pi_n(Y) \ \ |
\ \ [\alpha,
f_\#(\beta)]_w = 0 \hbox{ \ \ for all \ \
}\beta \in \pi_*(X)
\}.$$
Then $G_n(Y,X;f)$ is a subgroup of
$P_n(Y,X;f)$, and $P_*(X)$
occurs as the special case $X=Y$ and $f=1$.
In \thmref{thm:G
versus P}, we identify the quotient $P_*(Y,
X;f) / G_*(Y, X;f)$
for general maps $f\colon X \to Y$ between
rational spaces via a
particular commutative diagram of adjoints
and derivation spaces.
On the other hand, in \thmref{thm:coformal
G=P} we prove $G_*(Y,
X;f) = P_*(Y, X;f)$ when $f\colon X \to Y$
is a coformal map of
rational spaces.

As a second application, we describe and
study the rationalization
of the so-called \emph{$G$-sequence of a
map} as constructed by
Lee and Woo \cite{L-W1}. This is a chain
complex of the form
$$\xymatrix{\cdots \ar[r] & G_{n}(X)
\ar[r]^-{f_\#}& G_{n}(Y, X;f)
\ar[r] & G_{n}^{rel}(Y,X;f) \ar[r]&  G_{n-
1}(X) \ar[r]& \cdots}
$$
which arises naturally in the context of
function spaces and
evaluation maps described above. We identify
the rationalization
of the $G$-sequence in the framework of
derivation spaces and
adjoint maps of Quillen models
(\corref{cor:G-sequence}).  We
give, in particular, a complete calculation
of the homology of
this sequence at a particular term, in the
case of a single
cell-attachment (\thmref{thm:G-Sequence one
cell}). We also prove
the rationalized $G$-sequence is exact at
the same particular term
for coformal maps (\thmref{thm:G-sequence
exact coformal map}), as
an extension of the equality above.

The paper is organized as follows.
\secref{section:algebra}
contains purely algebraic definitions and
constructions used
throughout the paper. Here, we describe the
basic framework of
generalized derivation spaces and
generalized adjoint maps, and
the exact sequences in homology which arise
in this context.
In \secref{section:Quillen model} we obtain a  first  connection
between generalized
derivations and topology.
In a self-contained argument, we describe  an explicit Quillen minimal
model for a product $A
\times X$ in terms of the Quillen models of
the factors, in the
case in which $A$ is a rational
co-$H$-space.
The description depends on a particular class of generalized
derivations.
While we believe the general  result
is of independent interest,  we require it
here mainly for the case
in which $A$ is a single sphere $S^n$ and
briefly for the case in
which $A = S^n \vee S^n$. In
\secref{section:main result} we
connect the algebraic framework of \secref{section:algebra} to the
homotopy data of function spaces and evaluation maps. The main
results here are \thmref{thm:isos Phi and Phi*} and
\thmref{theorem:evaluation sequence}, in
which we identify the
rational homotopy groups of function spaces,
and more generally
the long exact rational homotopy sequence of
the evaluation
fibration $\omega\colon \map(X,Y;f) \to Y$,
in terms of homology
of derivation spaces of Quillen models.
Sections
\ref{section:evaluation group} and
\ref{section:G-sequence}
contain our examples and applications. In a
technical appendix, we
review some basic material concerning DG Lie
algebra homotopy
theory and provide details of some results
from this area that we
use for our proofs.

Throughout this paper, spaces $X$ and $Y$
will be based, simply
connected CW complexes of finite type. Given
a map $f \colon A \to
B$, either topological or algebraic, $f^*$
denotes pre-composition
by $f$ and $f_*$ post-composition by $f$.
We use ``$*$" to denote
the constant map of spaces. We write
$\omega$ generically for an
evaluation map $\omega \colon \map(X,Y;f)
\to Y$, since it will be
clear from the context which component is
intended.   We use
$H(f)$ to denote the map induced on homology
by $f$, and, when $f$
is a map of spaces, $f_\#$ for the map
induced on homotopy groups.
We write $f\sim g$ to denote that based maps
$f$ and $g$ are
homotopic via a based homotopy. In some
instances, we will
consider homotopies that are either unbased--- such as those in the
unbased function space component
$\map(X,Y;f)$, or are relative to
some subspace. In such instances, we will
specify the nature of
the homotopy.  We will continue to introduce
notation and set
conventions throughout the paper, on an as-needed basis and
usually at the start of each section.

\section{Lie Derivation Spaces}\label{section:algebra}

In this section, we consider a class of
chain complexes over $\Q$
obtained from generalized derivations of
differential graded Lie
algebras. Geometric content will be obtained
when we apply these
basic constructions to Quillen models in
\secref{section:main result}. We first establish our notational
conventions for this
material.

By \emph{vector space}, we mean a graded
vector space of finite
type over the rationals. Furthermore, all
our vector spaces will
be \emph{connected}, that is, they will be
positively graded. The
degree of an element $x \in V$ is written
$|x|$. The space of
degree $n$ elements of $V$ will be denoted
$V_n$.   A vector space
generated by a single element $v$ will be
written $\Q v$.  That
is, if $|v| = n$, then $(\Q v)_i$ is
isomorphic to $\Q$ if $i = n$
and is zero otherwise. The $k$th
\emph{suspension} of $V$, denoted
by $s^k(V)$, is the vector space defined as
$(s^k(V))_n =
V_{n-k}$. In particular, the
\emph{desuspension} of $V$ denoted
$s^{-1}V$ is the vector space with
$(s^{-1}V)_n = V_{n+1}$.  A
typical example of desuspension that we
consider is that of the
reduced rational homology of a simply
connected space of finite
type, denoted by $s^{-1}\widetilde{H}_*(X;\Q)$.  A vector space
$V$ equipped with a \emph{differential},
that is, a linear map
$d\colon V \to V$ of degree $-1$ that
satisfies $d^2 = 0$, will be
called  a \emph{DG (differential graded)
vector space}. Given a
cycle $x  \in V_n,$ a DG vector space,  we
write $\langle x\rangle$ for the homology class in $H_n(V, d)$
represented by $x$.

A \emph{Lie algebra} $L$ is a (graded)
vector space $L$ equipped
with a bilinear, degree zero map (the
``bracket" multiplication)
$[ \ \ , \ \ ] \colon L \times L \to L$ that
satisfies
\begin{enumerate}
\item[(i)] Anti-symmetry: $[\alpha, \beta] =
- (-1)^{|\alpha||\beta|}[\beta,
\alpha]$
\item[(ii)] Jacobi identity: $[\alpha,
[\beta, \gamma]] = [[\alpha,\beta], \gamma]
+
(-1)^{|\alpha||\beta|}[\beta, [\alpha,
\gamma]]$.
\end{enumerate}
The motivating example here is $\pi_*(\Omega X) \otimes \Q$, the
rational homotopy of $\Omega X$ with the
Samelson bracket $[ \ \ , \ \ ]$.

A \emph{DG Lie algebra} is a pair $(L, d)$
where $L$ is a Lie
algebra and $d$ is a vector space
differential that satisfies the
derivation law $d([x, y]) = [d(x), y] +(-1)^{|x|}[x, d(y)]$.  A
map of DG Lie algebras $\phi \colon (L, d)
\to (L', d')$ is a
homomorphism that respects differentials,
that is, satisfies $\phi
d = d' \phi$. Given a map of DG Lie
algebras, we may pass to
homology in the usual way.  If the induced
homomorphism $H(\phi)$
is an isomorphism, then we say $\phi$ is a
\emph{quasi-isomorphism}.  We will
frequently use the symbol
``$\simeq$" to denote the fact that a map is
a quasi-isomorphism,
especially in diagrams.  We write $\L(V)$
for the free Lie algebra
generated by the vector space $V$.  The
coproduct (or ``free
product") of (DG) Lie algebras $L$ and $L'$
is written $L\sqcup L'$.  We usually abuse notation somewhat and
write $\L(V,W)$ for
the coproduct $\L(V)\sqcup\L(W) = \L(V\oplus W)$. Similarly, we
will also write $\L(v)$ for $\L(\Q v)$,
$\L(v,w)$ for $\L(\Q v, \Q w) = \L (\Q v \oplus \Q w) =
\L(\Q v) \sqcup \L( \Q w)$,  $\L(v,W)$
for $\L(\Q v, W) = \L(\Q v \oplus W) = \L(\Q v) \sqcup \L(W)$, and
so-forth.  We commit another abuse of
notation by saying that a DG
Lie algebra is free if the underlying Lie
algebra is a free Lie
algebra.  To reduce parentheses, we usually
write a free DG Lie
algebra $\big(\L(V), d\big)$ as $\L(V;d)$.
Finally, a DG Lie
algebra $(L, d)$ is \emph{minimal} if $L$ is
free and the
differential is decomposable, that is, $d(L)
\subseteq [L, L]$.

We assume the reader is familiar with the
basic facts of rational
homotopy theory from the Quillen point of
view, that is, using DG
Lie algebra minimal models. Good references
for this material
include \cite{Tan83} and \cite[Part IV]{F-H-T}. Specifically, we
recall that each space $X$ has a
\emph{Quillen minimal model}
which is a minimal DG Lie algebra
$(\mathcal{L}_X, d_X)$ whose
isomorphism type is a complete invariant of
the rational homotopy
type of $X$. As a Lie algebra, we have
$\mathcal{L}_X =\L(s^{-1}\widetilde{H}_*(X;\Q))$.  The
differential $d_X$ is
determined by the topology of $X$ in a more
arcane way. The
Quillen model of $X$ recovers the rational
homotopy Lie algebra of
$X$ via a Lie algebra isomorphism
$H_*(\mathcal{L}_X, d_X) \cong
\pi_*(\Omega X) \otimes \Q$. Furthermore, a
map of spaces
$f\colon X \to Y$ induces a map of Quillen
models $\mathcal{L}_{f}\colon \mathcal{L}_X \to \mathcal{L}_Y$ with
$H(\mathcal{L}_f)$
corresponding, via the above isomorphism,
to the rationalization
of the map induced on homotopy groups by
$\Omega f\colon \Omega X
\to \Omega Y$.   We refer to this map of
Quillen models as the
\emph{Quillen minimal model} of the map $f$.
Our basic reference
for rational homotopy theory is \cite{F-H-T}.

Let $(L, d)$ be a DG Lie algebra. A
\emph{derivation of degree
$n$} of $L$ is a linear map
$\theta\colon L \to L$ that raises
degree by $n$ and satisfies the rule
$$ \theta([\alpha, \beta]) =
[\theta(\alpha), \beta] +
(-1)^{|\theta||\alpha|}[\alpha,
\theta(\beta)].$$
When $n =1$ we also require $d\circ\theta =
- \theta\circ d$. We
write $\Der_n(L)$ for the space of
derivations of degree $n$ of
$L$. The space $\Der_*(L)$ has the structure
of a DG Lie algebra,
with (commutator) bracket $[\theta_1,
\theta_2] = \theta_1 \circ
\theta_2 -(-1)^{|\theta_1||\theta_2|}
\theta_2 \circ \theta_1$,
and differential defined as $D(\theta) = [d,
\theta]$. The
\emph{adjoint map} $\ad\colon L \to
\Der_*(L)$, given by
$\ad(x)(y) = [x, y]$, is now a map of DG Lie
algebras.

At the expense of the Lie bracket, we can
extend the notion of a
derivation of a DG Lie algebra to a
derivation with respect to a
map of Lie algebras, as follows. Given a map
$\psi\colon (L, d_L)
\rightarrow (K,d_K)$ of DG Lie algebras,
define a \emph{derivation of degree $n$ with respect to $\psi$}, or
simply a
\emph{$\psi$-derivation of degree $n$}, to
be a linear map
$\theta\colon L \to K$ that increases degree
by $n$ and satisfies
$$\theta([\alpha, \beta]) = [\theta(\alpha),
\psi(\beta)] +
(-1)^{n|\alpha|}[\psi(\alpha),
\theta(\beta)]$$
for $\alpha, \beta \in L$.   Let $\Der_n(L,
K; \psi)$ denote the
space of all $\psi$-derivations of degree
$n$ from $L$ to $K$.
Next, define $D: \Der_n(L, K;\psi) \to
\Der_{n-1}(L, K;\psi)$ by
$D(\theta) = d_K \circ \theta  - (-1)^{|\theta|} \theta \circ
d_L$. The pair $\big(\Der_*(L, K;\psi), D\big)$ is then a DG
vector space.  The \emph{adjoint map associated to $\psi$} is
$\ad_{\psi}\colon K \to \Der_*(L, K; \psi)$
where
$$\ad_\psi(\alpha)(\beta) = [\alpha, \psi(\beta)].$$
It is easy to check that $\ad_\psi$ is a map
of DG vector spaces.
We continue to write $\ad$ for the adjoint
map associated to the
identity $1\colon L \to L$.

We now begin to consider homology of
derivation spaces.  At this
point, we introduce only what is needed for
\secref{section:main result}.  Further notions concerning
derivation spaces and their
homology will be introduced as needed in
subsequent sections. To
ease notation in what follows, we adopt the
following conventions.
If $(V, d)$ is a DG vector space we suppress
the differential when
we write its $n$th homology group, writing
just $H_n(V)$.
Furthermore, we indicate only the outermost
degree.  For example,
$H_n(\Der(L, K;\psi))$ denotes the $n$th
homology group of the DG
vector space $(\Der_*(L, K;\psi), D)$.

Next, we construct a long exact homology
sequence that we will
show corresponds to the long exact rational
homotopy sequence of
the evaluation fibration
(\thmref{theorem:evaluation sequence}).
Given a map $\phi\colon V \to W$ of DG
vector spaces we will need
the \emph{relativization of $\phi$}. This is
the DG vector space
$(\Rel_*(\phi), \delta)$ given, in degree
$n$, by
$$\Rel_n(\phi) = W_{n} \oplus V_{n-1}$$
with differential of degree $-1$ defined as
$$\delta(w, v) = (\phi(v) - d_W(w),
d_V(v)).$$
The inclusion $J: W_n \to \Rel_n(\phi)$ with
$J(w_n) = (w_n,0)$,
and the projection $P: \Rel_n(\phi) \to
V_{n-1}$ with
$P(w_n,v_{n-1}) = v_{n-1}$, give a short
exact sequence of chain
complexes
$$\xymatrix{0 \ar[r] & W_* \ar[r]^-{J} &
\Rel_*(\phi) \ar[r]^-{P} & V_{*-1} \ar[r]&
0}.$$
This leads to a long exact sequence whose
connecting homomorphism
is $H(\phi)$.  Actually, $J$ is really an
``anti-chain" map, in
that $\delta J = - J d_W$.  But $J$ still
induces a homomorphism
$H(J)$ on homology, and gives rise to a long
exact sequence in the
usual way. Thus we have a long exact
homology sequence
$$ \cdots  \to H_{n+1}(\Rel(\phi))
\stackrel{H(P)}{\Longrightarrow}  H_n(V)
\stackrel{H(\phi)}{\Longrightarrow}
H_n(W) \stackrel{H(J)}{\Longrightarrow}
H_n(\Rel(\phi)) \to
\cdots.$$
We refer to this sequence as the \emph{long
exact homology
sequence of $\phi$.}

Now apply this construction to the adjoint
map $ad_\psi\colon K
\to \Der_*(L, K;\psi)$ from above.  We
obtain a long exact
homology sequence
\begin{displaymath}
\xymatrix{ & \cdots \ar[r]^-{H(J)} &
H_{n+1}(\Rel(\ad_\psi))
 \ar `d[l]  `[lld]_(0.7){H(P)} [lld]
\\
H_{n}(K) \ar[r]^-{H(\ad_\psi)} &
H_{n}\big(\Der(L, K;\psi) \big)
\ar[r]^-{H(J)} & H_n(\Rel(\ad_\psi))
\ar`d[l] `[lld]_(0.7){H(P)}[lld] \\
H_{n-1}(K) \ar[r]^-{H(\ad_\psi)} &  \cdots
 }
\end{displaymath}
We call this sequence the \emph{long exact
derivation homology
sequence of $\psi$.}

\section{Quillen models for certain
products}\label{section:Quillen model}

The Sullivan minimal model of a product of
spaces is easily
expressed as the tensor product of the
Sullivan minimal models of
the spaces.  For Quillen minimal models, the
situation is more
complicated due to the fact that the direct
sum of minimal DG Lie
algebras---which is the categorical product
in this setting---is
not minimal. Tanr\'{e} has described how to
construct the Quillen
minimal model of a product of spaces $A
\times X$ in terms of the
Quillen minimal models of $A$ and $X$ (see
\cite[Prop.VII.1(2)]{Tan83}).  For our
purposes, we need explicit
minimal models for the cases in which $A =
S^n$ and $A = S^n\vee
S^n$.  In this section, we describe a
Quillen minimal model of $A
\times X$ for the more general case in which
$A$ is any rational
co-$H$-space. Our description essentially
makes Tanr{\'e}'s
construction explicit for these cases,
although our treatment here
is self-contained.

Suppose $X$ has Quillen minimal model
$\L(W;d_X)$, and that
$\L(V;d = 0)$ is the Quillen minimal model
of a simply connected,
finite-type rational co-$H$-space $A$.
Suppose $\{ v_i \}_{i \in
J}$ is a (connected, finite type) basis for
$V$, and $|v_i| = n_i
- 1$ for each $i \in J$. Topologically, this
corresponds to $A$
being of the rational homotopy type of the
wedge of spheres
$\bigvee_{i\in J} S^{n_i}$. We construct a
new minimal DG Lie
algebra from these data as follows: For each
$i \in J$, let $W_i$
denote the $n_i$-fold suspension of $W$,
that is, set $W_i =
s^{n_i}(W)$. Let $\lambda\colon \L(W)
\rightarrow \L(W, V,\oplus_{i\in J}W_i)$ denote the inclusion of
graded Lie algebras.
For each $i\in J$, let $S_i \colon W \to
W_i$ denote the
$n_i$-fold suspension isomorphism and define
a $\lambda$-derivation $S_i \colon \L(W) \to
\L(W, V, \oplus_{i\in J}W_i)$ by extending $S_i$ using the
$\lambda$-derivation rule.
Now define a differential $\partial$ on
$\L(W, V, \oplus_{i\in J}W_i)$ that extends the differentials on
$\L(W)$ and $\L(V)$,  by
setting
$$\partial(S_i(w)) = [v_i, w] + (-1)^nS_i(d_X(w))$$
for each generator $w \in W$ (and thus each
generator of $W_i$).
This definition is equivalent to the
identity $\ad(v_i) = D(S_i)$
in $\Der_{n_i-1}(\L(W), \L(W, V, \oplus_{i\in J}W_i);\lambda)$.

We will show that $\L(W, V, \oplus_{i\in J}W_i;\partial)$ is the
Quillen minimal model of $A \times X$. First
we check that the
preceding formula defines a differential.

\begin{lemma}
The derivation $\partial$ of $\L(W, V, \oplus_{i\in J}W_i)$
satisfies $\partial \circ\partial = 0$.
\end{lemma}

\begin{proof} It is sufficient to check on
generators.  Further, since
$\partial$ extends the differentials on
$\L(W)$ and $\L(V)$, it is
sufficient to check that
$(\partial)^2(S_i(w)) = 0$ for each $w
\in W$ and $i \in J$.  We compute as
follows:
$$
\begin{aligned}
(\partial)^2(S_i(w))  & =
\partial \left( [v_i, w] + (-1)^{n_i}
S_i(d_X(w)) \right)   \\
& = (-1)^{n_i-1}[v, d_X(w)]  + (-1)^{n_i}
\partial S_i(d_X(w)) \\
& = (-1)^{n_i-1}\ad(v_i)(d_X(w)) + (-1)^{n_i}\partial S_i(d_X(w))\\
& = S_id_X(d_X(w)) = 0.
\end{aligned}
$$
The penultimate step follows from the
identity $\ad(v_i) = D(S_i)$
observed above, which expands to give $[v_i, \chi] = \partial
S_i(\chi) - (-1)^{n_i} S_id_X(\chi)$ for any
$\chi \in \L(W)$, and
which here is applied to $d_X(w) \in \L(W)$.
\end{proof}

We will need the following technical point
in our argument.  The
discussion here is lifted from
\cite[Sec.22(f)]{F-H-T}. Suppose
that $\L(V;d)$ is a free, but not
necessarily minimal, connected
DG Lie algebra. That is, suppose that the
differential $d$ may
have a non-trivial linear part.  Since
$\L(V)$ is free, we can
write $d$ as the sum of two derivations $d =
d_0 + d_{+}$, where
$d_0\colon V \to V$ is the linear part of
the differential $d$ and
$d_{+}$ is the decomposable part of $d$ that
increases bracket
length. Indeed, $d_0$ itself is a
differential---although $d_{+}$
is generally not---and so we obtain a DG
vector space $(V,
d_{0})$. Next, suppose $\phi \colon \L(V;d')
\to \L(W;d)$ is a
morphism of free, but not necessarily
minimal, connected DG Lie
algebras. Again, because the Lie algebras
are free, we may write
the linear map $\phi\colon V \to \L(W)$ as a
sum $\phi = \phi_0 +
\phi_{+}$, where $\phi_0 \colon V \to W$ is
the linear part of
$\phi$ and $\phi_{+}$ is the decomposable
part of $\phi$ that
increases bracket length.  In this way, we
obtain a morphism of DG
vector spaces
$$\phi_0 \colon (V, (d')_{0}) \to (W, d_{0}),$$
which we refer to as the
\emph{linearization} of $\phi$.  The
following result can be interpreted as a
version of Whitehead's
theorem in our context.

\begin{lemma}\label{lem:linearization}
Let $\phi\colon \L(V;d') \to \L(W;d)$ be a
morphism of connected
free DG Lie algebras. Let $\phi_0 \colon
(V,(d')_0) \to (W, d_0)$
be the linearization of $\phi$.  Then $\phi$
is a
quasi-isomorphism of DG Lie algebras if, and
only if, $\phi_0$ is
a quasi-isomorphism of DG vector spaces.
\end{lemma}

\begin{proof} This is proved as
\cite[Prop.22.12]{F-H-T}.
\end{proof}

We now come to the main point of the
section.

\begin{theorem}\label{thm:Quillen model of product}
Let $X$ be a simply connected space of
finite type with Quillen
minimal model $\L(W;d_X)$.  Let $A$ be a
rational co-$H$-space of
the rational homotopy type of the wedge of
spheres $\bigvee_{i\in J} S^{n_i}$. Then
$\L(W, V, \oplus_{i\in J}W_i;\partial)$, as
described above, is the Quillen minimal
model of $A \times X$.
\end{theorem}

\begin{proof} Our starting point is the
well-known fact that the
direct sum of Quillen minimal models gives a
(non-minimal) DG Lie
algebra model for the product
\cite[p.332,Ex.3]{F-H-T}.  In our
case, this gives $\L(V)\oplus\L(W;d_X)$ as a
non-minimal model for
$A \times X$.  We will show that the obvious
projection
$$p\colon \L(W, V, \oplus_{i\in J}W_i;\partial) \to
\L(V)\oplus\L(W;\partial)$$
is a quasi-isomorphism.  Since the domain is
a minimal DG Lie
algebra, this is sufficient to show that it
is the Quillen minimal
model of the product.

So consider the following commutative
diagram of DG Lie algebra
morphisms:
$$\xymatrix{0 \ar[r] & K \ar[r]^-{i}\ar[d]_{p'} &
\L(W, V, \oplus_{i\in J}W_i;\partial)
\ar[r]^-{q}\ar[d]_{p} & \L(V) \ar[r]
\ar[d]_{1} & 0\\
0 \ar[r] & \L(W;d_X) \ar[r]^-{i'} &
\L(V)\oplus\L(W; d_X)\ar[r]^-{q'} & \L(V) \ar[r] & 0} $$
Here, $q$ and $q'$ are the obvious
(quotient) projections onto
$\L(V)$, and $i$ and $i'$ are the inclusions
of the kernels, so
that the rows are short exact sequences of
DG Lie algebras.  We
will argue that $p' \colon K \to \L(W;d_X)$
is a
quasi-isomorphism. First note that, as a
sub-DG Lie algebra of a
connected, free DG Lie algebra, $K$ is
itself a connected, free DG
Lie algebra. Indeed, as a Lie algebra, we
may write
$$K = \L(W, \oplus_{i}W_i, [V, W],
\oplus_{i}[V, W_i], [V,[V, W]],
\oplus_{i}[V,[V,W_i]],\dots;\partial_K)$$
or more succinctly $K =
\L(\{\ad^j(V)(W)\}_{j\geq0},
\{\oplus_{i}\ad^j(V)(W_i)\}_{j\geq0};\partial_K)$. In these
expressions, $[V,W]$ denotes the vector
space spanned by brackets
$[v,w]$ with $v \in V$ and $w \in W$, and
so-forth, and
$\ad^0(V)(W)$ denotes $W$, $\ad^1(V)(W)$
denotes $[V,W]$,
$\ad^2(V)(W)$ denotes $[V, [V,W]]$, and so-
forth.  We now claim
that $(\partial_K)_0$, the linear part of
the differential in $K$,
induces isomorphisms
$$(\partial_K)_0\colon
\oplus_{i\in J}\ad^j(V)(W_i) \to
\ad^{j+1}(V)(W)$$
for each $j\geq 0$.  First recall that
$\partial_K$ is simply the
restriction of the differential $\partial$
to the kernel of $q$,
and that $\partial(V) = 0$.   Extending our
notation a little
further, we can denote a typical spanning
element of $\ad^j(V)(W)$
as follows.  Suppose
$(v_{r_1},v_{r_2},\ldots, v_{r_j}) \in V^j$
is a $j$-tuple.  Then write
$\ad(v_{r_1},v_{r_2},\ldots,
v_{r_j})(w)$ for $[v_{r_1}, [v_{r_2},
[\ldots,[v_{r_{j-1}},v_{r_j}]] \ldots]$.  Likewise for elements
of $\ad^j(V)(W_i)$. Now
let $w \in W$ be a typical element.  From
the definition of
$\partial$ above, we have
$$
\begin{aligned}
\partial\big(\ad(v_{r_1},v_{r_2},\ldots,
v_{r_j})(S_i(w))\big) &= \pm
\ad(v_{r_1},v_{r_2},\ldots,
v_{r_j})(\partial(S_i(w))) \\
&= \pm \ad(v_{r_1},v_{r_2},\ldots, v_{r_j}, v_i)(w)\\
&\ \ \ \ \ \ \pm \ad(v_{r_1},v_{r_2},\ldots, v_{r_j})(S_id_X(w)).
\end{aligned}
$$
Since $\L(W;d_X)$ is the Quillen minimal
model of $X$, $d_X(w)$ is
decomposable in $\L(W)$ and thus $S_id_X(w)$
is decomposable in
$\L(W,W_i)$. From this, it follows that the
last term displayed
above, namely $\ad(v_{r_1},v_{r_2},\ldots,
v_{r_j})(S_id_X(w))$,
is decomposable in $K$. We prove this
assertion in
\lemref{lem:decomposable in K} below.
Assuming for the time being
its validity, it follows that the linear
part of the differential
in $K$ induces isomorphisms $(\partial_K)_0
\colon \ad^j(V)(W_i)
\cong \ad^{j}(V)\ad(v_i)(W)$ for each $i \in
J$ and each $j\geq
0$, and hence isomorphisms $(\partial_K)_0
\colon
\oplus_{i}\ad^j(V)(W_i) \cong
\ad^{j+1}(V)(W)$ for each $j\geq0$,
as claimed. Notice that as a consequence of
this, we must have
$(\partial_K)_0 = 0$ on each vector space of
generators
$\ad^{j+1}(V)(W)$ in $K$, for $j \geq 0$,
since the linear part of
a differential is itself a differential.  In
any case, this latter
fact also follows from
\lemref{lem:decomposable in K}. Finally,
notice that $(\partial_K)_0 = 0$ on the
vector space of generators
$W$ in $K$, since $\partial= d_X$ is
decomposable on $W$. It now
follows that the DG vector space $(Q(K),
(\partial_K)_0)$ obtained
by linearizing $K$ may be written as a
direct sum
$$(Q(K), (\partial_K)_0) \cong
(W,(\partial_K)_0 = 0)
\oplus_{j\geq0}\Big(\big(\oplus_{i}\ad^j(V)(W_i)\big) \oplus
\big(\ad^{j+1}(V)(W)\big),(\partial_K)_0 \Big),$$
in which each summand
$\Big(\big(\oplus_{i}\ad^j(V)(W_i)\big) \oplus
\big(\ad^{j+1}(V)(W)\big),(\partial_K)_0 \Big)$ is an
acyclic DG vector space.  It is now evident
that $H(Q(K),
(\partial_K)_0) \cong W$ and that the
linearization of $p'$, that
is, $(p')_0 \colon (Q(K), (\partial_K)_0)
\to (W,\partial_0 = 0)$,
is a quasi-isomorphism of DG vector spaces.
It follows from
\lemref{lem:linearization} that $p'$ is a
quasi-isomorphism of DG
Lie algebras.

Returning to the diagram of short exact
sequences, we now have
left and right vertical arrows that are
quasi-isomorphisms.
Therefore, by passing to the induced diagram
of long exact
homology sequences and applying the five-lemma, we obtain that $p$
is a quasi-isomorphism.  Hence $\L(W, V, \oplus_{i\in J}W_i;
\partial)$ is the Quillen minimal model of
$A \times X$.
\end{proof}

The proof of \thmref{thm:Quillen model of
product} will be
completed when we establish the following
lemma:

\begin{lemma}\label{lem:decomposable in K}
With notation as in the proof above, suppose
$\chi$ is a
decomposable element in $K$.  Then $[v, \chi]$ is also
decomposable in $K$, for any $v \in V$.  In
particular, if $\chi$
is decomposable in $\L(W, W_i)$ for some $i
\in J$, and hence
decomposable in $K$, then
$\ad(v_{r_1},v_{r_2},\ldots,v_{r_j})(\chi)$ is decomposable in $K$ for
any $j$-tuple
$(v_{r_1},v_{r_2},\ldots, v_{r_j}) \in V^j$.
\end{lemma}

\begin{proof}
Recall that $K = \L(\{\ad^j(V)(W)\}_{j\geq0}, \oplus_{i}
\{\ad^j(V)(W_i)\}_{j\geq0})$ is a sub-Lie
algebra of $\L(W,
V,\oplus_{i} W_i)$, and observe that $v \in
V$ is not a
generator---is not even an element---of $K$,
so the statement is
not entirely trivial. Without loss of
generality, we may assume
that $\chi$ is a monomial term.  We argue by
induction on the
bracket length in $K$ of $\chi$.  When
$\chi$ has length $2$ in
$K$, we have $\chi = [\chi_1, \chi_2]$ with
$\chi_1$ and $\chi_2$
indecomposable monomials in $K$.   From the
Jacobi identity in
$\L(W, V, \oplus_{i}W_i)$, we may write
$[v,\chi] = \pm [\chi_1,
[v, \chi_2]] \pm [\chi_2, [v, \chi_1]]$.
Now $\chi_1$ is an
element from either $\ad^j(V)(W)$ or
$\ad^j(V)(W_i)$ for some $i
\in J$ and $j\geq0$, and likewise for
$\chi_2$. So $[v, \chi_1]$
is from $\ad^{j+1}(V)(W)\subseteq K$ or
$\ad^{j+1}(V)(W_i)\subseteq K$, and likewise
for $[v, \chi_2]$.
That is, if $\chi_1$ and $\chi_2$ are
indecomposable in $K$, then
$[v,\chi] = \pm [\chi_1, [v, \chi_2]] \pm
[\chi_2, [v, \chi_1]]$
displays $[v, \chi]$ as decomposable (and
again of length $2$) in
$K$. Now assume inductively that the
assertion is true for $\chi$
of bracket length $\leq r$ in $K$. Let
$\chi$ be a monomial of
bracket length $r+1$ in $K$.  By judicious
use of the Jacobi
identity in $K$, we may assume that $\chi =
[\chi_1, \chi_2]$ with
$\chi_1$ an indecomposable in $K$ and
$\chi_2$ a decomposable
monomial of bracket length $r$ in $K$. Once
again, the Jacobi
identity in $\L(W, V, \oplus_{i}W_i)$ yields
$[v,\chi] = \pm
[\chi_1, [v, \chi_2]] \pm [\chi_2, [v, \chi_1]]$. Our inductive
hypothesis implies that $[v, \chi_2]$ is
(decomposable) in $K$,
and the same observations as were used to
start the induction show
that $[v, \chi_1]$ is an (indecomposable)
element in $K$.  Hence
$[v, \chi]$ is decomposable in $K$, and the
induction is complete.
The result follows.
\end{proof}

Since it is the main case we require here,
we write out explicitly
what this gives for the model of $S^n \times X$, with a slight
easing of notation.

\begin{corollary}\label{cor:Quillen model S times X}
Suppose $X$ has Quillen minimal model
$\L(W;d_X)$, let $\L(v)$
with $|v| = n-1$ and zero differential be
the Quillen model of
$S^n$, and set $W' = s^n(W)$.  Let $\lambda
\colon\L(W) \to \L(W,v, W')$ be the inclusion, and $S \colon\L(W)
\to  \L(W, v, W')$ be
the $\lambda$-derivation that extends the
linear map $S(w) = w'$.
Define a differential $\partial$ on $\L(W, v, W')$ by $\partial(w)
= d_X(w)$, $\partial(v) =0$, and
$$\partial(w') =
 [v, w] + (-1)^nS(d_X(w)),$$
for each $w \in W$.  Then $\L(W, v, W';\partial)$ is the Quillen
minimal model of $S^n \times X$.
\end{corollary}

\section{Lie  Derivations and Homotopy Groups of Function Spaces}%
\label{section:main result}

Say two maps of vector spaces $f\colon U \to
V$ and $g\colon U'
\to V'$ are \emph{equivalent} if there exist
isomorphisms $\alpha$
and $\beta$ which make the diagram
\begin{displaymath} \label{eq:partial square}
\xymatrix{U \ar[r]^{f}
\ar[d]_{\alpha}^{\cong} & V
\ar[d]_{\beta}^{\cong}\\
U' \ar[r]_{g} & V'}
\end{displaymath}
commutative.  We will extend this notion of
equivalence in the
obvious way to exact sequences of vector
spaces, and any other
diagram of vector space maps.  Given any map
$f\colon X \to Y$, we
have the homomorphism
$$j_\#\otimes1 \colon \pi_{n}\big(\map_*(X, Y;f)\big)\otimes\Q \to
\pi_{n}\big(\map(X,Y;f) \big)\otimes\Q$$
induced on rational homotopy groups by the
fibre inclusion of the
general evaluation fibration $\xymatrix{
\map_*(X, Y;f) \ar[r]^{j}
& \map(X, Y;f) \ar[r]^-{\omega} & Y}$.  On
the other hand, we have
the homomorphism
$$H(J) \colon
H_{n}\big(\Der(\mathcal{L}_X,\mathcal{L}_Y
;\mathcal{L}_f)\big)
\to  H_{n}(\Rel(\ad_{\mathcal{L}_f})\big)$$
that forms part of the long exact homology
sequence of the adjoint
map $\ad_{\mathcal{L}_f} \colon
\mathcal{L}_Y \to
\Der(\mathcal{L}_X,\mathcal{L}_Y ;\mathcal{L}_f)$ of the Quillen
minimal model $\mathcal{L}_f \colon
\mathcal{L}_X \to
\mathcal{L}_Y$ of $f$. In \thmref{thm:isos Phi and Phi*}, which is
the basic result of this section, we
establish that these two
homomorphisms are equivalent.  This result
and one immediate
consequence will occupy the remainder of
this section.

The main step is to establish vector space
isomorphisms
$$\Phi\colon \pi_n(\map_*(X,Y;f))\otimes\Q
\to H_n(\Der(\mathcal{L}_X,\mathcal{L}_Y;
 \mathcal{L}_f))$$
and
$$\Psi\colon \pi_n(\map(X,Y;f))\otimes\Q \to
H_n(\Rel(\ad_{\mathcal{L}_f})),$$
that give the equivalence. In the following,
we assume a fixed
choice of Quillen minimal model
$\mathcal{L}_f\colon \mathcal{L}_X \to \mathcal{L}_Y$ for $f\colon X \to Y$.
Write $\mathcal{L}_X =
\L(W;d_X)$  and $\mathcal{L}_Y = \L(V;d_Y)$.
To this end, we
define group homomorphisms
$$\Phi'\colon \pi_n(\map_*(X,Y;f)) \to
H_n(\Der(\mathcal{L}_X,\mathcal{L}_Y;
 \mathcal{L}_f))$$
and
$$\Psi'\colon \pi_n(\map(X,Y;f)) \to
H_n(\Rel(\ad_{\mathcal{L}_f}))$$
from the ordinary homotopy groups to the
appropriate vector
spaces, for $n\geq2$. Then the isomorphisms
$\Phi$ and $\Psi$ are
obtained as the rationalizations of these
homomorphisms.

Define $\Phi'$ as follows. Let $\alpha \in
\pi_n(\map_*(X,Y;f))$
be represented by a map $a\colon S^n \to
\map_*(X,Y;f)$.  Then the
adjoint $A \colon S^n \times X \to Y$ of $a$
has Quillen minimal
model $\mathcal{L}_{A}\colon
\mathcal{L}_{S^n \times X} \to
\mathcal{L}_Y$.  Recall from
\corref{cor:Quillen model S times X}
that $\mathcal{L}_{S^n \times X} = \L(W, v, W';\partial)$.  Now,
since $a$ is a (based) map into the function
space of \emph{based}
maps, we have $A\circ i_1 = *\colon S^n \to
Y$ and $A\circ i_2 =
f\colon X \to Y$. It follows that we may
take the Quillen minimal
model of $A$ to be a DG Lie algebra map
$$\mathcal{L}_{A}\colon  \L(W, v,W';\partial) \to \mathcal{L}_Y$$
that satisfies $\mathcal{L}_{A}(v) = 0$ and
$\mathcal{L}_{A}(w) =
\mathcal{L}_{f}(w)$ for each $w \in W$.  See
the appendix for
justification of this last assertion. Now
define a linear map
$\theta_A\colon \L(W) \to \mathcal{L}_Y$
that increases degree by
$n$ as the composition
$$\xymatrix{\L(W) \ar[r]^-{S} & \L(W, v, W')
\ar[r]^-{\mathcal{L}_{A}}&
\mathcal{L}_Y},$$
where $S\colon \L(W) \to \L(W, v, W')$ is
the derivation from
\corref{cor:Quillen model S times X}. A
straightforward check
shows that $\theta_A$ is an $\mathcal{L}_f$-derivation in
$\Der_n(\mathcal{L}_X,\mathcal{L}_Y;
\mathcal{L}_f)$. Furthermore,
we have $d_Y \theta_A = (-1)^n\theta_A d_X$
and so $\theta_A$ is a
cycle in the derivation space. Finally, we
set $\Phi'(\alpha) =
\langle \theta_A\rangle \in
H_n(\Der(\mathcal{L}_X,\mathcal{L}_Y;
\mathcal{L}_f))$. We will establish that
$\Phi'$ is a well-defined
homomorphism, and that its rationalization
$\Phi$ is an
isomorphism, in \thmref{thm:isos Phi and Phi*} below.

We define $\Psi'$, and thus its
rationalization $\Psi$, in a
similar manner. Let $\alpha \in \pi_n(\map(X,Y;f))$ be represented
by a map $a\colon S^n \to \map(X,Y;f)$.
Then the adjoint $A
\colon S^n \times X \to Y$ of $a$ still
satisfies $A\circ i_2 =
f\colon X \to Y$, but the composition
$A\circ i_1\colon S^n \to Y$
may give a non-trivial element of
$\pi_n(Y)$. Correspondingly, the
Quillen minimal model of $A$ satisfies
$\mathcal{L}_{A}(w) =
\mathcal{L}_{f}(w)$ for each $w \in W$, but
$\mathcal{L}_{A}(v)
\in \mathcal{L}_Y$ is now some non-trivial
$d_Y$-cycle. As before,
setting $\theta_A = \mathcal{L}_{A}\circ
S\colon \L(W) \to
\mathcal{L}_Y$ defines an
$\mathcal{L}_f$-derivation in
$\Der_n(\mathcal{L}_X,\mathcal{L}_Y;
\mathcal{L}_f)$.  Recalling
the definition of
$\Rel_*(\ad_{\mathcal{L}_f})$ from
\secref{section:algebra}, we obtain an
element
$$(\theta_A, \mathcal{L}_{A}(v)) \in
\Rel_n(\ad_{\mathcal{L}_f}).$$
Let $\delta$ and $D$ denote the
differentials in
$\Rel_*(\ad_{\mathcal{L}_f})$ and
$\Der_n(\mathcal{L}_X,\mathcal{L}_Y;
\mathcal{L}_f)$ respectively.
Then we have $\delta(\theta_A, \mathcal{L}_{A}(v)) = \big(
\ad_{\mathcal{L}_f}(\mathcal{L}_{A}(v)) -
D(\theta_A),
d_Y\mathcal{L}_{A}(v))\big)$, and
$d_Y\mathcal{L}_{A}(v) = 0$.  We
check that
$\ad_{\mathcal{L}_f}(\mathcal{L}_{A}(v)) -
D(\theta_A)
= 0 \in \Der_{n-1}(\mathcal{L}_X,\mathcal{L}_Y;
\mathcal{L}_f)$.
For $\chi \in \mathcal{L}_X$, we have
$$
\begin{aligned}
\ad_{\mathcal{L}_f}(\mathcal{L}_{A}&(v))(\chi) - D(\theta_A)(\chi)  \\
&= [\mathcal{L}_A(v), \mathcal{L}_f(\chi)] -
d_Y\big(\mathcal{L}_{A}\circ S(\chi)\big) +
(-1)^{n}\mathcal{L}_{A}\circ Sd_X(\chi)  \\
&= \mathcal{L}_{A}\big(
[v,\chi] + (-1)^{n} Sd_X(\chi) \big) -
d_Y\mathcal{L}_{A}\big(S(\chi)\big)  \\
&= \mathcal{L}_{A}\partial\big(S(\chi)\big) -
d_Y\mathcal{L}_{A}\big(S(\chi)\big)\\
&=0.
\end{aligned}
$$
Thus $(\theta_A, \mathcal{L}_{A}(v))$ is a
cycle in the relative
chain complex. Finally, we set
$\Psi'(\alpha) = \langle \theta_A,
\mathcal{L}_{A}(v)\rangle \in
H_n\big(\Rel(\ad_{\mathcal{L}_f})\big)$.

In the following result, we establish the
basic properties of
$\Phi$ and $\Psi$.

\begin{theorem}\label{thm:isos Phi and Phi*}
Suppose $n \geq 2$. Then we have:
\begin{enumerate}
\item[(A)] $\Phi'$ and $\Psi'$ are well-
defined homomorphisms;
\item[(B)] Their rationalizations $\Phi$ and
$\Psi$ are
isomorphisms;
\item[(C)] The following square is
commutative:
$$\xymatrix{\pi_n(\map_*(X,Y;f))\otimes\Q
\ar[r]^-{\Phi}_-{\cong}\ar[d]_{j_\#\otimes\Q} &
H_n(\Der(\mathcal{L}_X,\mathcal{L}_Y;
\mathcal{L}_f))
\ar[d]_{H(J)}\\
\pi_n(\map(X,Y;f))\otimes\Q \ar[r]_-{\Psi}^-{\cong}&
H_n(\Rel(\ad_{\mathcal{L}_f}))}$$
\end{enumerate}
\end{theorem}

\begin{proof}
Throughout the proof we will give full
details for arguments
concerning $\Phi$.  The arguments for $\Psi$
are similiar, and we
will simply indicate them without details.
We will need to use
some facts about homotopy in the DG Lie
algebra setting. The most
complete reference for this material is
\cite[Ch.II.5]{Tan83}. The
appendix to this paper contains a quick
overview, and also
provides careful justifications of some
technical details used in
the following proof.  We will make free use
of the notation
concerning DG Lie algebra homotopy reviewed
in the appendix.

(A) \emph{$\Phi'$ is well-defined}:
Suppose that $a\sim b\colon
S^n \to \map_*(X,Y;f)$ are two
representatives of the homotopy
class $\alpha$.  The adjoint of the (based)
homotopy in
$\map_*(X,Y;f)$ from $a$ to $b$ gives a
homotopy of their
adjoints, $A, B\colon S^n \times X \to Y$,
that is stationary on
the subset $S^n \vee X \subseteq S^n \times
X$. Indeed, the
homotopy of the adjoints is stationary at
the constant map on
$S^n$ and is stationary at $f$ on $X$.
Consequently, the
corresponding Quillen minimal models
$\mathcal{L}_{A},
\mathcal{L}_{B}\colon \mathcal{L}_{S^n \times X} \to
\mathcal{L}_Y$ are homotopic via a DG Lie
algebra homotopy
$$\mathcal{H} \colon \L(W, v,W')_I \to \mathcal{L}_Y$$
that satisfies $\mathcal{H}(\L(v)_I) = 0$,
$\mathcal{H}(\L(sW,
\widehat{W})) = 0$, and $\mathcal{H}(w) =
\mathcal{L}_{f}(w)$ for
each $w \in W$ (see \lemref{lem:restricted DG Lie homotopy} for
details). Let $\sigma \colon \L(W, v, W')_I
\to \L(W, v, W')_I$ be
the derivation of degree $1$ defined in the
appendix.  Define a
linear map $\Theta \colon \L(W) \to
\mathcal{L}_Y$ of degree $n+1$
as the composition $\Theta =
\mathcal{H}\circ\sigma\circ S$.  A
straightforward check, using the fact that
$\mathcal{H}(sW,
\widehat{W}) = 0$, shows that $\Theta$ is an
$\mathcal{L}_f$-derivation.  We will show
that $D\Theta = \theta_B
- \theta_A \in \Der_n(\mathcal{L}_X,
\mathcal{L}_Y;
\mathcal{L}_f)$.  We have $\theta_B =
\mathcal{H}\circ
\lambda_1\circ S$ and $\theta_A =
\mathcal{H}\circ \lambda_0\circ
S = \mathcal{H}\circ S$. Let $J$ denote the
ideal of $\L(W, v,
W')_I$ generated by $sW\oplus\widehat{W}$.
Since $\sigma$ vanishes
on the generators of $J$, and $D_I$
preserves the set of
generators, $J$ is stable under the
composition $\sigma D_I$.
Furthermore, $\mathcal{H}$ is zero on $J$,
since it is zero on the
generators.  We claim that $(\sigma D_I)^r
\colon W' \to \L(W, v,W')_I$ has image in $J$, for each $r \geq 2$.
First observe that
$$\begin{aligned}
\sigma D_I(w') = \sigma \partial (w') &=
\sigma ( [v,w] + (-1)^n S
d_X(w) )\\
& = [sv,w] + (-1)^n [v, sw] + (-1)^n \sigma
Sd_X(w).
\end{aligned}$$
Furthermore, the only terms not in $J$ that
may appear in $\sigma
S d_X(w)$ are terms in the sub-Lie algebra
$\L(W, sW')$ that have
exactly one occurrence of an element from
$sW'$.  On applying
$D_I$ to such terms, the only terms still
not in $J$ that may
appear in $D_I \sigma  S d_X(w)$ are terms
in the sub-Lie algebra
$\L(W, sW', \widehat{W'})$ that have exactly
one occurrence of an
element either from $sW'$, or from
$\widehat{W'}$.  Since $\sigma
$ is zero on $sW'\oplus\widehat{W'}$,
$\sigma D_I\sigma S d_X(w)
\in J$. Direct computation shows that both
$\sigma D_I([sv,w])$
and $\sigma D_I([v, sw])$ are in $J$.   That
is, $(\sigma
D_I)^r(w') \in J$, for each $r \geq 2$.
Consequently,
$\mathcal{H}$ vanishes on these terms.
Finally, using this, we
compute that
$$\begin{aligned}
\theta_B(w) &= \mathcal{L}_B\circ S(w) =\mathcal{H}\circ\lambda_1(w')\\
&= \mathcal{H}\big(w' + D_I\sigma (w') +
\sigma D_I(w') +
\sum_{r\geq2}
\frac{1}{r!} (\sigma D_I)^r (w')\big)\\
&= \mathcal{H}(w')  + \mathcal{H}D_I\sigma
(w') + \mathcal{H}
\big([sv,w] + (-1)^n [v, sw] + (-1)^n \sigma
Sd_X(w)\big)\\
&= \mathcal{L}_A(w') + d_Y \mathcal{H}\sigma
(w') + (-1)^n
\mathcal{H}\sigma S d_X(w) \\
&= \theta_A(w) + d_Y \Theta(w) - (-1)^{n+1}
\Theta d_X(w).
\end{aligned}
$$
It follows that the difference of
derivations $\theta_B - \theta_A
= D\Theta$ in $\Der_n(\mathcal{L}_{X},
\mathcal{L}_{Y};
\mathcal{L}_f)$. Hence $\Phi'$ is well-
defined.

\emph{$\Phi'$ is a homomorphism}: Let $\nu
\colon S^n \to S^n \vee
S^n$ denote the usual pinching
comultiplication. Given $\alpha,
\beta \in \pi_{n}(\map_*(X,Y;f))$,  the sum
$\alpha+ \beta$ is the
composition $(\alpha\mid\beta)\circ\nu$.
Suppose $\alpha, \beta$
have adjoints $A, B \colon S^n \times X \to
Y$, respectively.  Let
$i_1, i_2 \colon S^n \to S^n\vee S^n$ denote
the inclusions, and
let $(A\mid B)_f\colon (S^n \vee S^n)\times
X \to Y$ be the map
defined by $(A\mid B)_f \circ(i_1\times 1) =
A$ and $(A\mid B)_f
\circ(i_2\times 1) = B$. Then the adjoint of
$\alpha+ \beta$ is $C
:= (A\mid B)_f\circ(\nu\times1) \colon
S^n\times X \to Y$. We
focus on identifying the Quillen minimal
model of $(A\mid B)_f$,
and it will follow that $\Phi'_f$ is a
homomorphism.

The map $(A\mid B)_f\colon (S^n \vee
S^n)\times X \to Y$ is
determined, up to homotopy, as the unique
map $F$ that makes the
following diagram homotopy commutative:
$$\xymatrix{ S^n\times X \ar[d]_{i_1\times1}
\ar[rrd]^{A}\\
(S^n\vee S^n)\times X \ar@{.>}[rr]^-{F}& &
Y\\
S^n\times X \ar[u]^{i_2\times1} \ar[rru]_{B}
}
$$
Consequently, any DG Lie algebra map
$\Gamma$ that makes the
diagram
$$\xymatrix{\L(W, v,
W';\partial)\ar[d]_{j_1}\ar[rrd]^{\mathcal{L}_A}\\
\L(W, v_1, v_2, W_1, W_2;
\partial)\ar@{.>}[rr]_(0.65){\Gamma} &
&\mathcal{L}_Y \\
\L(W, v, W';\partial) \ar[u]^{j_2}
\ar[rru]_{\mathcal{L}_B} }
$$
commute up to DG homotopy is a Quillen model
for $(A\mid B)_f$. In
this diagram, $\L(W, v_1, v_2, W_1, W_2;
\partial)$ is the Quillen
model of $(S^n \vee S^n)\times X$ as
described in
\thmref{thm:Quillen model of product}, and
$j_1, j_2$ the obvious
inclusions that identify $v$ with $v_1$ and
$W'$ with $W_1$, and
$v$ with $v_2$ and $W'$ with $W_2$
respectively.   This
characterization of (the Quillen model of)
$(A\mid B)_f$ is
explained in detail in \cite{L-S}. Now there
is an obvious choice
for $\Gamma\colon \L(W, v_1, v_2, W_1, W_2;
\partial) \to \mathcal{L}_Y$, namely the map
that makes the
diagram commute.  Hence, this map is a
Quillen model for $(A\mid
B)_f$. Finally, since $\nu\colon S^n \to S^n
\vee S^n$ has Quillen
model $\mathcal{L}_{\nu} \colon \L(v) \to
\L(v_1, v_2)$ given by
$\mathcal{L}_{\nu}(v) = v_1 + v_2$, it
follows that $\nu\times
1\colon S^n\times X \to (S^n \vee S^n)\times
X$ has Quillen model
$\mathcal{L}_{\nu\times 1} \colon \L(W, v, W';\partial)
\to \L(W,v_1, v_2, W_1, W_2; \partial)$ given by
$\mathcal{L}_{\nu\times1}(w') = w_1 + w_2$.
Thus
$$\mathcal{L}_{C}(w') =
\Gamma\circ\mathcal{L}_{\nu\times1}(w') =
\Gamma(w_1 + w_2) = \theta_A(w)
+ \theta_B(w).$$
Hence we have $\theta_C =
\mathcal{L}_{C}\circ S = \theta_A +
\theta_B$ and it follows that $\Phi'$ is a
homomorphism.

\emph{$\Psi'$ is a well-defined
homomorphism}:  This is
established by making small adjustments to
the preceding arguments
for $\Phi'$.  In this case, the homotopy of
the adjoints $A$ and
$B$ is stationary at $f$ on $X$, but is not
stationary on $S^n$.
The corresponding homotopy $\mathcal{H}$ of
Quillen minimal models
still satisfies $\mathcal{H}(\L(sW,
\widehat{W})) = 0$, and
$\mathcal{H}(w) = \mathcal{L}_{f}(w)$ for
each $w \in W$ (see
\lemref{lem:restricted DG Lie homotopy} for
details), but we must
allow for non-zero terms in
$\mathcal{H}(\L(v)_I)$. Since $v$ is a
$\partial$-cycle, however, we have
$\lambda_1(v) = v + \widehat{v}
\in \L(W, v, W')_I$.  Therefore,
$\mathcal{L}_B(v) =
\mathcal{H}\circ\lambda_1(v) =
\mathcal{H}(v) +
\mathcal{H}D\sigma(v) = \mathcal{L}_A(v) +
d_Y(\mathcal{H}(sv))$.
Since $\mathcal{H}$ still vanishes on
$sW\oplus\widehat{W}$, the
composition $\Theta =
\mathcal{H}\circ\sigma\circ S$ still defines
an $\mathcal{L}_f$-derivation. Computing
exactly as before, we
find that $\theta_B(w) = \theta_A(w) +
D\Theta (w) +
[\mathcal{H}(sv),\mathcal{L}_f(w)]$. It
follows that, in the
relative complex, we have $\delta(-
\Theta,\mathcal{H}(sv)) =
(\theta_B - \theta_A, \mathcal{L}_B(v) -
\mathcal{L}_A(v))$.  So
$\Psi'$ is well-defined. To show that
$\Psi'$ is a homomorphism,
we can use precisely the same argument as
for $\Phi'$, and use the
additional identity
$$\mathcal{L}_{C}(v) =
\Gamma\circ\mathcal{L}_{\nu\times1}(v) =
\Gamma(v_1 + v_2)
= \mathcal{L}_A(v) + \mathcal{L}_B(v)$$
at the final step.

(B) \emph{$\Phi$ is a surjection}: Suppose
given $\theta \in
\Der_n(\mathcal{L}_X,\mathcal{L}_Y;
\mathcal{L}_f)$, a cycle
derivation of degree $n$.  Define a Lie
algebra map
$\mathcal{L}_{A}\colon  \L(W, v, W';
\partial) \to \mathcal{L}_Y$
by setting
$$\mathcal{L}_A(w) =
\mathcal{L}_f(w), \ \ \ \mathcal{L}_A(v) = 0 \ \  \hbox{\ \ and \ \ }
 \mathcal{L}_{A}(w') = \theta(w)$$
for $w \in W$.  Just as in the definition of
$\Phi'$,
$\mathcal{L}_A \circ S$ is an
$\mathcal{L}_f$-derivation and by
construction we have $\mathcal{L}_A \circ S
= \theta$. We check
that $\mathcal{L}_A$ commutes with
differentials as follows:
$$ \begin{aligned}
\mathcal{L}_A(\partial(w')) & =
\mathcal{L}_A([v, w] +
(-1)^{n}S(d_{X}(w)) \\
&= 0 + (-1)^{n} \mathcal{L}_A(S(d_{X}(w)) \\
& =   (-1)^{n} \theta(d_{X}(w)) \\
& = d_Y(\theta(w)) \\
&= d_Y ( \mathcal{L}_A(w'))
\end{aligned}
$$
Let $A\colon S^n \times X \to Y_\Q$ be the
geometric realization
of $\mathcal{L}_A$, from the correspondence
between (homotopy
classes of) maps between rational spaces and
DG Lie algebra maps
between Quillen models.  Let  $i_1\colon S^n
\to S^n \times X$ and
$i_2\colon  X \to S^n \times X$ denote the
inclusions. Since
$\mathcal{L}_A \circ \mathcal{L}_{i_1} = 0$
and $\mathcal{L}_A
\circ \mathcal{L}_{i_2} = \mathcal{L}_f$, we
have $A \circ i_1
\sim *$ and $A \circ i_2 \sim f_\Q$.
Altering the geometric
realization $A$ up to homotopy, we may
assume $A \circ i_1 = *$
and $A \circ i_2 = f_\Q$. Thus, the adjoint
$a\colon S^n \to
\map_*(X,Y_\Q;f_\Q)$ of $A$ represents an
element $\alpha \in
\pi_n(\map_*(X,Y_\Q;f_\Q))$. Clearly, we
have $\Phi(\alpha) =
\langle \theta \rangle$, and so $\Phi$ is
surjective.

\emph{$\Phi$ is an injection}: Since $\Phi$
is a homomorphism, it
is sufficient to check that $\Phi(\alpha) =
0$ implies $\alpha = 0
\in \pi_n(\map_*(X,Y_\Q;f_\Q))$. As before,
write $\Phi(\alpha) =
\langle \theta_A\rangle$ and suppose
$\theta_A \in
\Der_n(\mathcal{L}_X, \mathcal{L}_Y;
\mathcal{L}_f)$ is a boundary
so that $\theta_A = D(\Theta)$ for $\Theta
\in
\Der_{n+1}(\mathcal{L}_X, \mathcal{L}_Y;
\mathcal{L}_f)$. Define a
homotopy
$$\mathcal{G}\colon  \L(W, v, W')_I \to
\mathcal{L}_Y$$
by setting $\mathcal{G} = \mathcal{L}_A$ on
$\L(W, v, W')$ (so
$\mathcal{G}$ starts at $\mathcal{L}_A$),
$\mathcal{G} = 0$ on
$\L(sW, \widehat{W}, sv, \widehat{v})$,  and
$\mathcal{G}(sw') =
\Theta(w)$ on generators $w' \in W'$. We
then set
$\mathcal{G}(\widehat{w'}) = \mathcal{G}(D_I sw') = d_Y
\mathcal{G}(sw') = d_Y\Theta(w)$, for
$\widehat{w'} \in
\widehat{W'}$, and then extend $\mathcal{G}$
as a Lie algebra map,
so that $\mathcal{G}$ is a DG Lie algebra
map.   It is
straightforward to check that $\mathcal{G}$
ends at
$\mathcal{G}\circ\lambda_1(w) =
\mathcal{L}_f(w)$ and
$\mathcal{G}\circ\lambda_1(v) = 0$ on
$\L(W,v)$.  Just as in the
proof of (A) above, we observe that because
$\mathcal{G}$ is zero
on $\L(sW, \widehat{W})$, the composition
$\mathcal{G}\circ \sigma
\circ S$ acts as a derivation in
$\Der_{n+1}(\mathcal{L}_X,
\mathcal{L}_Y; \mathcal{L}_f)$. Therefore,
we have
$\mathcal{G}\circ \sigma \circ S (d_X w) =
\Theta(d_Xw)$ for each
$w \in W$.  Furthermore, again just as in
part (A),  we find that
$\mathcal{G}$ is zero on terms $(\sigma
D_I)^r(w')$, for each $r
\geq 2$. Therefore, for $w' \in W'$, this
homotopy ends at
$$\begin{aligned}
\mathcal{G}\circ\lambda_1(w') &=
\mathcal{G}\big(w' + D_I \sigma
(w') + \sigma D_I(w') + \sum_{r\geq2}
\frac{1}{r!} (\sigma D_I)^r (w')\big)\\
&= \mathcal{G}(w')  + \mathcal{G}D_I \sigma
(w') + \mathcal{G}
\sigma((-1)^n
S d_X(w) )\\
&= \mathcal{G}(w') + d_Y \mathcal{G}\sigma
(w') + (-1)^n
\mathcal{G} \sigma S d_X(w) \\
&= \theta_A(w) + d_Y \Theta(w) - (-1)^{n+1}
\Theta d_X(w),
\end{aligned}
$$
which is zero. Hence $\mathcal{G}$ ends at
the Quillen minimal
model of the composition $f \circ p_2 \colon
S^n \times X \to Y$.
It follows that $A \sim f \circ p_2 \colon
S^n \times X \to Y_\Q$.
Taking adjoints, we obtain that $a \sim *
\colon S^n \to
\map_*(X,Y;f)$. Actually, we only obtain
this last homotopy as a
free homotopy by taking adjoints, since the
homotopy between $A$
and $f\circ p_2$ is based, but not
necessarily relative to $X$.
However, a based map from a sphere is based-homotopic to the
constant map if it is freely homotopic to
the constant map
\cite[p.27]{Spa89}. Thus $\Phi$ is
injective.

\emph{$\Psi$ is an isomorphism}: Once again,
we need only make
slight adjustments to the preceding
arguments for $\Phi$. Suppose
given $(\theta, y)$ a $\delta$-cycle in
$\Rel_n(\ad_{\mathcal{L}_f})$.  Then $d_Y(y)
= 0$ and $D(\theta) =
\ad_{\mathcal{L}_f}(y)$, that is, we have
$$d_Y\theta(\chi) - (-1)^n \theta d_X(\chi)
= [y,\mathcal{L}_f(\chi)]$$
for $\chi \in \mathcal{L}_X$. In this case,
define
$\mathcal{L}_{A}\colon \L(W, v, W'; \partial) \to \mathcal{L}_Y$
by $\mathcal{L}_A(w) = \mathcal{L}_f(w)$,
$\mathcal{L}_A(v) = y$,
and $\mathcal{L}_{A}(w') = \theta(w)$.  Now
argue that $\Psi$ is
surjective following the same steps as were
taken for $\Phi$.

To show injectivity of $\Psi$, suppose that
$(\theta_A,\mathcal{L}_A(v)) \in
\Rel_n(\ad_{\mathcal{L}_f})$ is a boundary
in the relative complex, so that $(\theta_A,
\mathcal{L}_A(v)) =
\delta(\Theta, y)$.  That is, $\theta_A =
\ad_{\mathcal{L}_f}(y) -
d_Y\Theta + (-1)^{n+1}\Theta d_X$ and
$d_Y(y) = \mathcal{L}_A(v)$.
In this case, define the DG homotopy
$\mathcal{G}$ by setting
$\mathcal{G} = \mathcal{L}_A$ on $\L(W, v, W')$, $\mathcal{G} = 0$
on $\L(sW, \widehat{W})$, $\mathcal{G}(sv) = -y$,
$\mathcal{G}(\widehat{v}) = -
\mathcal{L}_A(v)$, $\mathcal{G}(sw')
= \Theta(w)$, and $\mathcal{G}(\widehat{w'})
= d_Y\Theta(w)$.
With this definition,  the same steps as
were used for $\Phi$ may
now be followed to show that $\Psi$ is
injective.

(C) The commutativity of the diagram follows
directly from the
definitions. Notice that $j \colon
\map_*(X,Y;f) \to
\map_*(X,Y;f)$ is the fibre inclusion of the
evaluation fibration.
Therefore, if $a$ is a representative of
$\alpha \in
\pi_n(\map_*(X,Y;f))$, we may also take $a$
to be a representative
of $j_\#(\alpha)\in \pi_n(\map(X,Y;f))$.
Consequently, the Quillen
minimal model of the adjoint of a
representative of $j_\#(\alpha)$
may be may be taken as $\mathcal{L}_A$ with
$\mathcal{L}_A(v) =
0$.  Hence, before rationalizing, we have
$\Psi'(j_\#(\alpha)) =
\langle \theta_A, 0\rangle = \langle
J(\theta_A)\rangle =
H(J)(\langle\theta_A\rangle) =
H(J)\circ\Phi'(\alpha)$.
\end{proof}

\begin{remark}\label{rem:naturality}
The homomorphisms $\Phi'$ and $\Psi'$, and
hence their
rationalizations, have certain naturality
properties.
Post-composition by a map $g \colon Y \to Z$
gives a map $g_*
\colon \map_*(X,Y;f) \to \map_*(X,Z;g\circ
f)$.  On the other
hand, the Quillen minimal model of $g$
induces a chain map
$(\mathcal{L}_g)_* \colon
\Der_*(\mathcal{L}_X, \mathcal{L}_Y;
\mathcal{L}_f) \to \Der_*(\mathcal{L}_X,
\mathcal{L}_Z;
\mathcal{L}_{g\circ f})$.  In this latter
chain complex, we choose
$\mathcal{L}_{g}\circ \mathcal{L}_{f}$ as
the Quillen minimal
model for $g\circ f$. If we denote the
homomorphism $\Phi$ of
\thmref{thm:isos Phi and Phi*} by
$\Phi^f_*$, so as to include the
original map in the notation, then we have
the following
commutative diagram:
$$\xymatrix{\pi_n\big(\map_*(X, Y;
f)\big)\otimes\Q \ar[r]^-
{\Phi^f_*}\ar[d]_{(g_*)_\#\otimes1} &
H_n\big(\Der(\mathcal{L}_{X},\mathcal{L}_{Y}
;\mathcal{L}_{f})\big)
\ar[d]^{H((\mathcal{L}_{g})_*)} \\
\pi_n\big(\map_*(X, Z; g\circ
f)\big)\otimes\Q
\ar[r]^-{\Phi^{g\circ f}_*} &
H_n\big(\Der(\mathcal{L}_X,
\mathcal{L}_Z; \mathcal{L}_{g\circ f})\big)}$$
Pre-composition by a map $h\colon W \to X$
gives a similar
naturality property of $\Phi$.  Furthermore,
$\Psi$ is natural in
the same sense.
\end{remark}

Schlessinger and Stasheff, and Tanr\'{e}
have constructed  a
(non-minimal) Quillen model for the
universal fibration
$$X \to B\,\map_*(X,X;1) \to
B\,\map(X,X;1)$$
in terms of Lie derivations and adjoint maps
(\cite{Sch-St,Tan83}).  Their result
specializes to identify the
long exact sequence induced in rational
homotopy groups by the
universal fibration, in the framework of
Quillen models and
derivations.  The following consequence of
\thmref{thm:isos Phi
and Phi*} extends this specialization of
their result, in that it
identifies the long exact rational homotopy
sequence of a general
evaluation fibration.

We first make the following observation.

\begin{lemma}\label{lemma:alpha prime}
Suppose given a diagram of vector spaces
$$\xymatrix{
B_{n+1} \ar[r]^{j_{n+1}}
\ar[d]_{\cong}^{\beta_{n+1}} &
C_{n+1}\ar[r]^{k_{n+1}}\ar[d]_{\cong}^{\gamma_{n+1}} & A_n
\ar[r]^{i_n}\ar@{.>}[d]^{\alpha_n} & B_n
\ar[r]^{j_n}
\ar[d]_{\cong}^{\beta_n}  &  C_n
 \ar[d]_{\cong}^{\gamma_n}   \\
Y_{n+1} \ar[r]_{q_{n+1}} &
Z_{n+1}\ar[r]_{r_{n+1}} & X_n
\ar[r]_{p_{n}} & Y_n \ar[r]_{q_{n}} &  Z_n
}$$
for each $n \geq 2$ (that is, a ``ladder
with every third rung
missing"). Suppose the rows are exact, each
$\beta_{n}$ and
$\gamma_{n}$ is an isomorphism, and
$\gamma_{n}\circ j_{n} =
q_{n}\circ\beta_{n}$ for each $n$. Then
there exist isomorphisms
$\alpha_{n}\colon A_n \to X_n$, for
$n\geq2$, which makes the
entire ladder commutative.
\end{lemma}

\begin{proof} This is straightforward.  See
\cite[Lem.3.1]{L-S} for details.
\end{proof}

Now consider the evaluation fibration
$\xymatrix@C=20pt{ \map_*(X,
Y;f) \ar[r]^{j} & \map(X, Y;f) \ar[r]^-
{\omega} & Y}$ for a map
$f\colon X \to Y$. By going one step back in
the Barratt-Puppe
sequence, we obtain a fibration
$$
\xymatrix{\Omega Y \ar[r]^-{\partial} &
\map_*(X,Y;f) \ar[r]^-{j}
&\map(X,Y;f)}
$$
and a long exact sequence in homotopy
\begin{equation}\label{eq:evaluation sequence}
\xymatrix@C=16pt{\cdots \ar[r]^-
{\partial_\#} &
\pi_{n}\left(\map_*(X, Y;f) \right)
\ar[r]^-{j_\#}  &
\pi_{n}\big(\map(X, Y;f) \big) \ar[r]^-
{\omega_\#} &
\pi_{n-1}(\Omega Y) \ar[r]^-{\partial_\#} &
\cdots}
\end{equation}
Here, we are identifying $\pi_{n-1}(\Omega
Y)$ with $\pi_{n}(Y)$
in the usual way. Notice that when $X = Y$
and $f=1$ the long
exact sequence (\ref{eq:evaluation
sequence}) is equivalent to
that of the universal fibration \cite{D-Z}.

\begin{theorem}\label{theorem:evaluation sequence}
Let $f\colon X \to Y$ be a map between
simply connected CW
complexes of finite type with $X$ finite.
Then the
rationalization of the long exact homotopy
sequence
(\ref{eq:evaluation sequence}), as far as
the term
$\pi_{1}\big(\Omega Y\big)\otimes\Q$, is
equivalent to the long
exact derivation homology sequence of the
Quillen minimal model
$\mathcal{L}_f\colon \mathcal{L}_X \to
\mathcal{L}_Y$ of $f$, that
is,
\begin{displaymath}
\xymatrix@C=20pt{\cdots \ar[r]^-
{H(\ad_{\mathcal{L}_f})} &
H_{n}\big(\Der(\mathcal{L}_X,\mathcal{L}_Y
;\mathcal{L}_f) \big)
 \ar[r]^-{H(J)} &
H_{n}(\Rel(\ad_{\mathcal{L}_f})\big)
\ar[r]^-{H(P)} &
H_{n-1}(\mathcal{L}_Y) \ar[r]^-
{H(\ad_{\mathcal{L}_f})} & \cdots   }
\end{displaymath}
as far as the term $H_{1}(\mathcal{L}_Y)$.
\end{theorem}

\begin{proof}
Replace the diagrams
$$\xymatrix{
B_n \ar[r]^{j_n} \ar[d]_{\cong}^{\beta_n}  &
C_n
 \ar[d]_{\cong}^{\gamma_n}   \\
Y_n \ar[r]_{q_n} &  Z_n }$$
of \lemref{lemma:alpha prime} with the
diagrams
$$\xymatrix{
\pi_n(\map_*(X,Y;f))\otimes\Q
\ar[r]^{j_\#\otimes\Q}
\ar[d]_{\cong}^{\Phi} &
\pi_n(\map(X,Y;f))\otimes\Q
 \ar[d]_{\cong}^{\Psi}   \\
H_{n}\big(\Der(\mathcal{L}_X,\mathcal{L}_Y
;\mathcal{L}_f) \big)
\ar[r]^{H(J)} &
H_{n}(\Rel(\ad_{\mathcal{L}_f})\big) }$$
for $n \geq 2$.  Then with the top row
replaced by the
rationalization of the long exact homotopy
sequence
(\ref{eq:evaluation sequence}), and the
bottom row by the long
exact derivation homology sequence of
$\mathcal{L}_f\colon
\mathcal{L}_X \to \mathcal{L}_Y$,
\lemref{lemma:alpha prime}
obtains isomorphisms
$$\alpha_n \colon \pi_n(\Omega Y)\otimes\Q
\to H_n(\mathcal{L}_Y)$$
with appropriate commutativity properties
for $n \geq 2$.  At the
bottom end, we have a diagram
$$\xymatrix{
\pi_2(\map_*(X,Y;f))\otimes\Q
\ar[r]^{j_\#\otimes\Q}
\ar[d]_{\cong}^{\Phi} &
\pi_2(\map(X,Y;f))\otimes\Q
 \ar[d]_{\cong}^{\Psi} \ar[r]^-
{\omega_\#\otimes\Q} &\pi_1(\Omega
Y)\otimes\Q \ar@{.>}[d]^{\alpha_1}   \\
H_{2}\big(\Der(\mathcal{L}_X,\mathcal{L}_Y
;\mathcal{L}_f) \big)
\ar[r]^{H(J)} &
H_{2}(\Rel(\ad_{\mathcal{L}_f})\big)
\ar[r]_{H(P)} & H_{1}(\mathcal{L}_Y) }$$
Here, we may use the standard identification
of $\pi_1(\Omega
Y)\otimes\Q$ with $H_{1}(\mathcal{L}_Y)$ for
our last isomorphism
$\alpha_1$.  Then we need only check that
the final square
commutes. For this, we regard the standard
identification
$\alpha_1 \colon \pi_1(\Omega Y)\otimes\Q
\to
H_{1}(\mathcal{L}_Y)$ as follows. Suppose $a
\colon S^2 \to Y$ is
a representative of $\eta \in \pi_1(\Omega
Y)\otimes\Q =
\pi_2(Y)\otimes\Q$.  Then we have
$\mathcal{L}_a \colon
\mathcal{L}_{S^2} \to \mathcal{L}_Y$ and we
set $\alpha_1(\eta) =
\mathcal{L}_a(u)$, where $\mathcal{L}_{S^2}
= \L(u)$ with $|u| =
1$.  Now suppose given $\zeta \in
\pi_2(\map(X,Y;f))\otimes\Q$
with adjoint $G \colon S^2 \times X \to
Y_\Q$.  Then
$\omega_\#(\zeta) = [G\circ i_1] \in
\pi_2(Y)\otimes\Q$.  Hence,
we have $H(P)\circ \Psi(\zeta) = H(P)
\langle \theta_G,
\mathcal{L}_G(v)\rangle = \langle
\mathcal{L}_G(v)\rangle =
\langle \mathcal{L}_G\circ
\mathcal{L}_{i_1}(u)\rangle =
\alpha\circ w_\#(\zeta)$. Thus we have the
equivalence asserted.
\end{proof}

\begin{remark}
The last part of the proof above extends to
give
$\alpha_n\circ(\omega_\#\otimes\Q) =
H(P)\circ\Psi$ for all $n
\geq 2$, if $\alpha_n$ is taken to be the
standard identification
$\alpha_n \colon \pi_n(\Omega Y)\otimes\Q
\to
H_{n}(\mathcal{L}_Y)$. However, we have
chosen to argue as above
to avoid the details involved in showing the
remaining squares of
the diagram commute with the standard choice
of $\alpha_n$.
Unfortunately, this means that we have not
established a
\emph{natural} equivalence of sequences in
\thmref{theorem:evaluation sequence}.
\end{remark}

\section{Evaluation Subgroups and Gottlieb's
Question}%
\label{section:evaluation group}

As an immediate consequence of
\thmref{theorem:evaluation
sequence}, we identify the rationalized
evaluation subgroup of a
map $f\colon X \to Y$ in terms of Quillen
models. Let $\psi \colon
K \to L$ be a DG Lie algebra map.

\begin{definition}\label{definition:evaluation group}
The \emph{$nth$ evaluation subgroup of the
DG  Lie algebra map
$\psi\colon L \to K$} is the subgroup
$$G_n(K,L; \psi) = \ker \{ H(\ad_\psi)\colon
H_{n-1}(K) \to H_{n-1}(\Der(L, K ;
\psi)) \}$$
of $H_{n-1}(K)$.  Notice the shift in
degrees.  We specialize this
to define the Gottlieb group of a DG Lie
algebra following
Tanr\'{e}.  The \emph{$nth$ Gottlieb group
of the DG  Lie algebra
$(L,d)$} is the subgroup
$$G_n(L) = \ker \{ H(\ad)\colon H_{n-1}(L)
\to H_{n-1}(\Der(L)) \}$$
of $H_{n-1}(L)$.
\end{definition}

\begin{theorem} \label{theorem:evaluation group}
Let $f\colon X \to Y$ be a map between
simply connected CW
complexes of finite type with $X$ finite.
Then the
rationalization of the evaluation subgroup
of $f$ is isomorphic to
the evaluation subgroup of the Quillen model
$\mathcal{L}_f\colon
\mathcal{L}_X \to \mathcal{L}_Y$ of $f$ as
in
\defref{definition:evaluation group}. That
is,  for $n \geq 2$ we
have
$$G_n(Y,X;f) \otimes \Q \cong   \ker\{
H(\ad_{\mathcal{L}_f})\colon H_{n-
1}(\mathcal{L}_Y)  \to H_{n-
1}(\Der(\mathcal{L}_X,
\mathcal{L}_Y; \mathcal{L}_f)) \}.$$
\end{theorem}

\begin{proof} From the long exact homotopy
sequence of the evaluation fibration,
the rationalized evaluation subgroup
$G_{n}(Y, X; f) \otimes \Q$
corresponds to the kernel of $\partial\colon
\pi_{n-1}(\Omega Y)
\otimes \Q \to \pi_{n-1}(\map_*(X, Y; f))
\otimes \Q$.  The result
follows from \thmref{theorem:evaluation
sequence}.
\end{proof}

Specializing to the identity map we recover

\begin{corollary}[{\cite[Cor.VII.4(10)]{Tan83}}]
    \label{cor:Tanre}
Let $X$ be a simply connected, finite
complex.  Then the
rationalization of the Gottlieb group of $X$
is the Gottlieb group
of the Quillen model $(\mathcal{L}_X, d_X)$
of $X$ as in
\defref{definition:evaluation group}.  That
is, for $n \geq 2$ we
have
$$G_n(X) \otimes \Q \cong   \ker\{
H(\ad_{\mathcal{L}_X})\colon H_{n-
1}(\mathcal{L}_X)  \to
H_{n-1}(\Der(\mathcal{L}_X)) \}.$$
\end{corollary}

We apply these identifications to address
Gottlieb's question on
the difference between the Whitehead
centralizer $P_*(X)$ and the
Gottlieb group $G_*(X)$.  In ordinary
homotopy theory,
constructing spaces with $G_*(X) \neq
P_*(X)$ represents a
challenging problem.  Ganea gave the first
example of inequality
in \cite{Gan}.  See \cite{Oprea} for a
recent reference and some
interesting examples of $G_1(X) \neq P_1(X)$
with $X$ a finite
complex.

Recall the definition of the generalized
Whitehead center of a map
$P_n(Y,X;f) \subseteq \pi_n(Y)$ from the
introduction. From the
identification of the Samelson product in
$\pi_*(\Omega Y)\otimes
\Q$ with the product in $H(\mathcal{L}_Y)$,
we have
$$P_{n}(Y_\Q, X_\Q;f_\Q) \cong \ker \{
 \ad_{H(\mathcal{L}_f)}\colon
H_{n-1}(\mathcal{L}_Y) \to
 \Der_{n-1}(H(\mathcal{L}_X),
H(\mathcal{L}_Y); H(\mathcal{L}_f))
 \}.$$
Notice that although the relative evaluation
subgroup behaves well
with respect to rationalization, in the
sense that $G_n(Y,
X;f)\otimes\Q = G_*(Y_\Q, X_\Q;f_\Q)$ (at
least for $X$ finite),
the inclusion $P_*(X) \otimes \Q \subseteq
P_*(X_\Q)$ is usually
strict.

Rationally, the difference between the $n$th
evaluation subgroup
of a map and the generalized Whitehead
center can be described
precisely. The difference is governed by the
``induced
derivation'' map
$$I\colon H_n(\Der(\mathcal{L}_X,
\mathcal{L}_Y;
\mathcal{L}_f)) \to \Der_n(H(\mathcal{L}_X),
H(\mathcal{L}_Y);
H(\mathcal{L}_f))$$
which we now introduce.

A $D$-cycle $\theta \in \Der_n(L, K; \psi)$
commutes (in the
graded sense) with the differentials of $L$
and $K$, and so
induces a map $H(\theta) \in \Der_n(H(L),
H(K); H(\psi))$ defined
by $H(\theta)(\langle \xi \rangle) = \langle
\theta(\xi) \rangle$
for $\xi$ a cycle of $L$. If $\theta$ is a
$D$-boundary then it
carries cycles of $L$ to boundaries of $K$.
Thus we obtain a
linear map
$$I\colon H_n(\Der(L, K; \psi)) \to
\Der_n(H(L),
H(K); H(\psi))$$
given by $I(\langle \theta \rangle) =
H(\theta)$ for $\theta$ a
cycle of $\Der_n(L, K; \psi)$.

Now consider the commutative diagram
$$ \xymatrix{
H_n(\mathcal{L}_Y) \ar[rr]^-
{H(\ad_{\mathcal{L}_f})}
\ar[drr]_-{\ad_{H(\mathcal{L}_f)}} & &
H_n(\Der(\mathcal{L}_X,
\mathcal{L}_Y;\mathcal{L}_f))
\ar[d]^{I} \\
& &  \Der_n(H(\mathcal{L}_X),
H(\mathcal{L}_Y);H(\mathcal{L}_f)).} $$

\begin{theorem}\label{thm:G versus P}
Let $f\colon X \to Y$ be a map between
simply connected CW
complexes of finite type with $X$ a finite
complex. For $n \geq
1$, we have
$$\frac{P_{n+1}(Y_\Q,
X_\Q;f_\Q)}{G_{n+1}(Y_\Q,X_\Q;f_\Q)} \cong
\ker(I) \cap \im( H(\ad_{\mathcal{L}_f})
).$$
\end{theorem}

\begin{proof}
The map $H(\ad_{\mathcal{L}_f})$ induces a
map
$$\overline{H(\ad_{\mathcal{L}_f})} \colon
\frac{\ker(
\ad_{H(\mathcal{L}_f)})}{\ker(H(\ad_{\mathcal{L}_f}))}
\to \ker(I) \cap \im( H(\ad_{\mathcal{L}_f})
)$$
that is easily checked to be an isomorphism.
The result follows
from \thmref{theorem:evaluation group} and
the above discussion,
in which we identify
$\ker(H(\ad_{\mathcal{L}_f}))$ with
$G_{n+1}(Y_\Q,X_\Q;f_\Q)$, and $\ker(
\ad_{H(\mathcal{L}_f)})$
with $P_{n+1}(Y_\Q, X_\Q;f_\Q)$.
\end{proof}

Using these notions, it is straightforward
to give examples of
maps $f\colon X \to Y$ which give an
inequality $G_*(Y_\Q, X_Q;
f_\Q) \neq P_*(Y_\Q, X_\Q; f_\Q)$.

\begin{example}
Suppose that $f \colon X \to Y$ is a
rationally trivial map.  Then
$\mathcal{L}_f = 0\colon \mathcal{L}_X \to
\mathcal{L}_Y$.  It
follows that $\ad_{\mathcal{L}_f} = 0\colon
\mathcal{L}_Y \to
\Der_*(\mathcal{L}_X,
\mathcal{L}_Y;\mathcal{L}_f)$.  In this
case, we have $G_*(Y_\Q, X_Q; f_\Q) =
P_*(Y_\Q, X_\Q; f_\Q) =
\pi_*(Y)\otimes\Q$.  However, suppose
that $f_\#\otimes\Q = 0 \colon
\pi_*(X)\otimes\Q \to
\pi_*(Y)\otimes\Q$, or even just that
$\im(f_\#\otimes\Q) \subseteq
P_*(Y_\Q)$. Then
$\ad_{H(\mathcal{L}_f)} = 0\colon
H_*(\mathcal{L}_Y) \to
\Der_*(H(\mathcal{L}_X),
H(\mathcal{L}_Y);H(\mathcal{L}_f))$ and so  $P_*(Y_\Q,
X_\Q; f_\Q) = \pi_*(Y)\otimes\Q$.  However,
the equality
$G_*(Y_\Q, X_Q; f_\Q) = \pi_*(Y)\otimes\Q$
would only hold if
$H(\ad_{\mathcal{L}_f}) = 0$, and there is
no particular reason
why this should be so.  To illustrate this
last case, consider
$f\colon \C P^2 \to S^4$ obtained by
pinching out the $2$-cell of
$\C P^2$.  This map has Quillen minimal
model
$$\mathcal{L}_f \colon \L(x_1, x_3;d_X) \to
\L(u_3;d_Y = 0)$$
with $\mathcal{L}_f(x_1) = 0$ and
$\mathcal{L}_f(x_3) = u_3$.
Here, the subscript of a generator denotes
its degree and
$d_X(x_1) = 0$, $d_X(x_3) = [x_1, x_1]$.
Now
$\ad_{\mathcal{L}_f}(u_3) \in
\Der_3(\mathcal{L}_X,
\mathcal{L}_Y;\mathcal{L}_f)$ is defined by
$\ad_{\mathcal{L}_f}(u_3)(x_1) = 0$ and
$\ad_{\mathcal{L}_f}(u_3)(x_3) = [u_3,
u_3]$.  On the other hand,
$D = 0$ in $\Der_3(\mathcal{L}_X,
\mathcal{L}_Y;\mathcal{L}_f)$.
Therefore, $H(\ad_{\mathcal{L}_f})(\langle
u_3) \rangle \not= 0$
and $G_4(Y_\Q, X_Q; f_\Q) = 0$.
\end{example}

We next give a class of examples for which
the equality $G_*(Y_\Q,
X_\Q;f_\Q) = P_*(Y_\Q, X_\Q;f_\Q)$ holds.
For this, we review the
notion of \emph{coformality} and some
terminology associated with
this concept. Suppose that a minimal DG Lie
algebra $\L(V;d)$ has
a second (or ``upper'') grading on the
generating subspace $V
=\oplus_{i\geq0} V^i$.  This extends to a
second grading of
$\L(V)$ in the obvious way, and we write
$\L(V)^i$ for the
sub-vector space of $\L(V)$ consisting of
all elements of $\L(V)$
of second grading equal to $i$. We also
write $V^{(i)}$ for the
sub-vector space of $V$ consisting of all
elements of $V$ of
second grading less than or equal to $i$.
Then we say that
$\L(V;d)$ is a \emph{bigraded} minimal DG
Lie algebra if the
differential decreases second degree
homogeneously by one, that
is, if $d(V^0) = 0$ and $d(V^i) \subseteq
\L(V)^{i-1}$ for $i \geq
1$.  If $\L(V;d)$ is a bigraded minimal DG
Lie algebra, then the
second grading passes to homology, making
$H(\L(V;d))$ a bigraded
Lie algebra.  We write $H^i(\L(V;d))$ for
the sub-vector space of
$H(\L(V;d))$ consisting of homology classes
represented by cycles
of upper degree equal to $i$, and we have
$H(\L(V;d)) =
\oplus_{i\geq0}H^i(\L(V;d))$.

\begin{definition}
    \label{definition:coformal Lie algebra}
Let $\L(V;d)$ be a bigraded minimal DG Lie
algebra in the above
sense.  We say $\L(V;d)$ is \emph{coformal}
if $H^i(\L(V;d)) = 0$
for $i > 0$, so that $H(\L(V;d)) =
H^0(\L(V;d))$.  We say that a
space $X$ is a \emph{coformal space} if its
Quillen minimal model
is coformal.
\end{definition}

Equivalently, we may define $\L(V;d)$ to be
coformal if there
exists a quasi-isomorphism $\rho\colon
\L(V;d) \to (H(L), d=0)$.
The equivalence of these definitions is
established by the notion
of a \emph{bigraded model} in the DG Lie
algebra setting---see
\cite[Sec.3]{H-S} for the Sullivan model
setting. There are many
interesting examples of coformal spaces:
Moore spaces and more
generally rational co-H-spaces, including
suspensions; some
homogeneous spaces; products and wedges of
coformal spaces.

This notion of coformality extends to a map.
Suppose that
$\phi\colon \L(V;d) \to \L(W;d')$ is a map
of bigraded minimal DG
Lie algebras as defined above.  If
$\phi(V^i) \subseteq \L(W)^i$
for each $i\geq0$, then we say that $\phi$
is a bigraded map.

\begin{definition}
    \label{definition:coformal map}
A map $\phi\colon \L(V;d) \to \L(W;d')$ of
bigraded minimal DG Lie
algebras is a \emph{coformal map} if both
$\L(V;d)$ and $\L(W;d')$
are coformal, and $\phi$ is a bigraded map
(with respect to the
second gradings that display the coformality
of $\L(V)$ and
$\L(W)$). A map of coformal spaces $f \colon
X \to Y$ is a
coformal map if its Quillen minimal model
$\mathcal{L}_f \colon
\mathcal{L}_X \to \mathcal{L}_Y$ is a
coformal map of bigraded
minimal DG Lie algebras.
\end{definition}

Equivalently, we may define $\phi\colon
\L(V;d) \to \L(W;d')$  to
be coformal if there exist quasi-isomorphisms $\rho\colon \L(V;d)
\to H(\L(V;d))$ and $\rho'\colon \L(W;d')
\to H(\L(W;d'))$ such
that the diagram
$$\xymatrix{\L(V;d) \ar[r]^{\phi}
\ar[d]_{\rho}^{\simeq} &
\L(W;d')\ar[d]^{\rho'}_{\simeq}\\
H(\L(V;d)) \ar[r]_{H(\phi)} & H(\L(W;d'))}$$
is DG homotopy commutative.

\begin{remark}
Suppose given a map of DG Lie algebras $\Phi
\colon L \to L'$.  By
constructing DG Lie algebra minimal models
of $L$ and $L'$, and
then using the standard lifting lemma, we
obtain a DG homotopy
commutative diagram
$$\xymatrix{\L(V;d) \ar[r]^{\phi}
\ar[d]_{\rho}^{\simeq} &
\L(W;d')\ar[d]^{\rho'}_{\simeq}\\
L \ar[r]_{\Phi = H(\phi)} & L'}.$$
In this way, it is always possible to
realize a Lie algebra map as
the homomorphism induced on homology by a
coformal map. Notice,
however, that there may be many other DG Lie
algebra maps $\L(V;d)
\to \L(W;d')$ that induce the same
homomorphism on homology as
$H(\phi)$.  Coformality, therefore,
distinguishes a unique DG
homotopy class of maps from amongst the
various realizations. From
this point of view, a coformal map is the
simplest realization of
its induced homomorphism on homology.
\end{remark}

\begin{theorem}\label{thm:coformal G=P}
Let $f\colon X \to Y$ be a coformal map
between CW complexes of
finite type with $X$ finite.  Then
$P_*(Y_\Q, X_\Q; f_\Q) =
G_*(Y_\Q, X_\Q; f_\Q)$.
\end{theorem}

\begin{proof} Suppose $\mathcal{L}_X =
\L(W;d_X)$ and
$\mathcal{L}_Y = \L(V;d_Y)$ are coformal,
and that $\mathcal{L}_f$
is bigraded.   Take $\alpha \in
H_n(\mathcal{L}_Y)$.   With
reference to \thmref{thm:G versus P}, we
show that if $I\circ
H(\ad_{\mathcal{L}_f})(\alpha) = 0$, then
$H(\ad_{\mathcal{L}_f})(\alpha) = 0$.  Since
$Y$ is coformal, we
may assume $\alpha = \langle \xi \rangle$
for a $d_Y$-cycle $\xi
\in \L(V)^0$. Observe that
$\ad_{\mathcal{L}_f}(\xi) \in
\Der_n(\mathcal{L}_X, \mathcal{L}_Y;
\mathcal{L}_f)$ is then a
$D$-cycle that preserves upper degree. If
$I\circ
H(\ad_{\mathcal{L}_f})(\alpha) = 0$, then
for each $d_X$-cycle
$\chi \in \L(W)$, we have
$\ad_{\mathcal{L}_f}(\xi)(\chi) =
d_Y(\eta)$ for some $\eta \in \L(V)$.  We
now use this to
construct $\theta \in
\Der_{n+1}(\mathcal{L}_X, \mathcal{L}_Y;
\mathcal{L}_f)$ such that $D(\theta) =
\ad_{\mathcal{L}_f}(\xi)$.

Since $X$ is coformal, we have $W =
\oplus_{i\geq0} W^i$, and each
$w \in W^0$ is a $d_X$-cycle.  Therefore, we
have
$\ad_{\mathcal{L}_f}(\xi)(w) = d_Y(\eta)$
for some $\eta \in
\L(V)$.  Furthermore, since $\mathcal{L}_Y$
is bigraded, we may
choose $\eta \in \L(V)^1$.  Use this to
define a linear map
$\theta_0 \colon W^0 \to \L(V)^1$ and extend
to an
$\mathcal{L}_f$-derivation $\theta_0 \in
\Der_{n+1}(\L(W^0),
\mathcal{L}_Y; \mathcal{L}_f)$. By
construction, we have
$D(\theta_0)(\chi) = d_Y(\theta_0(\chi)) =
\ad_{\mathcal{L}_f}(\xi)(\chi)$ for $\chi
\in \L(W^0)$.

Assume inductively that $\theta_m \in
\Der_{n+1}(\L(W^{(m)}),
\mathcal{L}_Y; \mathcal{L}_f)$ is defined,
increasing upper degree
homogeneously by $1$, and satisfying
$D(\theta_m) =
\ad_{\mathcal{L}_f}(\xi)$ on $\L(W^{(m)})$.
For $w \in W^{m+1}$,
consider the element
$\ad_{\mathcal{L}_f}(\xi)(w) + (-1)^{n+1}
\theta_m(d_X w)$.  Since
$\ad_{\mathcal{L}_f}(\xi)$ is a
$D$-cycle, and $d_Xw \in \L(W)^{m} \subseteq
\L(W^{(m)}$, we
compute that
$$\begin{aligned}
d_Y(\ad_{\mathcal{L}_f}(\xi)(w) + (-
1)^{n+1}& \theta_m(d_X w)) =
d_Y\ad_{\mathcal{L}_f}(\xi)(w) + (-1)^{n+1}
d_Y\theta_m(d_X w)\\
&= (-1)^n \ad_{\mathcal{L}_f}(\xi)(d_Xw) +
(-1)^{n+1}
\ad_{\mathcal{L}_f}(\xi)(d_X w) \\
&\ \ \ \ \ \ + \ad_{\mathcal{L}_f}(\xi)
d_X(d_X w)\\
& = 0.
\end{aligned}
$$
Since $\ad_{\mathcal{L}_f}(\xi)(w) + (-
1)^{n+1} \theta_m(d_X w)$
is a $d_Y$-cycle in $\L(V)^{m+1}$, we again
use the coformality of
$\mathcal{L}_Y$---specifically, that
$H^+(\mathcal{L}_Y) = 0$---to
conclude that there exists some $\zeta \in
\mathcal{L}_Y$ with
$d_Y(\zeta) = \ad_{\mathcal{L}_f}(\xi)(w) +
(-1)^{n+1}
\theta_m(d_X w)$. Furthermore, we may choose
$\zeta \in
\L(V)^{m+2}$. Clearly, $\zeta$ may be chosen
so as to depend
linearly on $w$.  So use this to define a
linear map $\theta_{m+1}
\colon W^{m+1} \to \L(V)^{m+2}$ with
$\theta_{m+1}(w) = \zeta$,
and extend $\theta_m$ to an $\mathcal{L}_f$-
derivation
$\theta_{m+1} \in \Der_{n+1}(\L(W^{(m+1)}),
\mathcal{L}_Y;
\mathcal{L}_f)$. By construction, we have
$$\begin{aligned}
D \theta_{m+1}(w)  & = d_Y\theta_{m+1}(w) -
(-1)^{n+1}
\theta_{m+1}(d_X w)\\
&= d_Y\zeta - (-1)^{n+1} \theta_{m}(d_X w) =
\ad_{\mathcal{L}_f}(\xi)(w) \end{aligned}
$$
for $w \in W^{m+1}$ and since
$D(\theta_{m+1})$ is an
$\mathcal{L}_f$-derivation, this gives
$D(\theta_{m+1})(\chi) =
\ad_{\mathcal{L}_f}(\xi)(\chi)$ for $\chi
\in \L(W^{(m+1)})$. This
completes the induction and gives an
$\mathcal{L}_f$-derivation
$\theta \in \Der_{n+1}(\L(W), \mathcal{L}_Y;
\mathcal{L}_f)$ that
satisfies $D(\theta) =
\ad_{\mathcal{L}_f}(\xi)$. The result
follows.
\end{proof}

The following special case is well-known.

\begin{corollary}\label{cor:coformal G=P}
Let $X$ be a simply connected, finite CW
complex. If $X$ is
coformal, then $P_*(X_\Q) = G_*(X_\Q)$.
\end{corollary}

\begin{proof}
Restrict \thmref{thm:coformal G=P} to the
case $f = 1\colon X \to
X$.
\end{proof}

\section{The Rationalized $G$-Sequence}\label{section:G-sequence}

In this section, we identify the
rationalization of the
$G$-sequence mentioned in the introduction.
Suppose given a map
$f \colon X \to Y$.  Then we have the
commutative square
\begin{equation}\label{eq:fibre connecting
maps}
\xymatrix{\Omega X \ar[r]^{\Omega f}
\ar[d]_{\partial}&  \Omega
Y\ar[d]^{\partial}\\
\map_*(X, X;1) \ar[r]^-{f_*} &
\map_*(X,Y;f),}
\end{equation}
in which the vertical maps are the
connecting maps arising from
the evaluation fibrations  $\omega\colon
\map(X,X;1) \to X$ and
$\omega: \map(X,Y;f) \to Y$ as in
\secref{section:main result}.
The maps $\Omega f$ and $f_*$ lead to long
exact homotopy
sequences and the vertical maps give
homomorphisms of
corresponding terms, yielding a homotopy
ladder
\begin{equation}\label{eq:homotopy ladder}
\xymatrix{\cdots \ar[r]^{p} & \pi_{n}(\Omega
X)
\ar[r]^{(\Omega f)_{\#}}
\ar[d]_{\partial_\#}& \pi_{n}(\Omega Y)
\ar[r]^{j}\ar[d]_{\partial_\#} &
\pi_n(\Omega f)
\ar[r]^{p}\ar[d]_{\partial_\#} & \cdots \\
\cdots  \ar[r]^-{\widehat{p}} &
\pi_{n}(\map_*(X,X;1))
\ar[r]^-{(f_*)_{\#}} &
\pi_{n}(\map_*(X,Y;f))
\ar[r]^-{\widehat{j}} &  \pi_n(f_*) \ar[r]^-
{\widehat{p}} &
\cdots}
\end{equation}
in the usual way. Whenever we have such a
ladder, with exact rows,
there is an associated ``kernel sequence,"
that is, a sequence
obtained by restricting the maps in the top
row to the kernels of
the vertical rungs. The $G$-sequence of the
map $f$ may be
defined, with a shift in degree, as the
kernel sequence of the
above homotopy ladder. A portion of this
construction is shown
here:
$$
\xymatrix@C=12pt{\cdots \ar[r]^-{p}&
   G_{n+1}(X) \ar@{^{(}->}[d]
\ar[r]^{(\Omega f)_{\#}}
  & G_{n+1}(Y, X;f) \ar@{^{(}->}[d]
\ar[r]^{j} & G_{n+1}^{rel}(Y, X;
  f) \ar@{^{(}->}[d] \ar[r]^-{p} & \cdots\\
  \cdots \ar[r]^-{p}
   & \pi_{n}(\Omega X)
\ar[d]^{\partial_{\#}}
  \ar[r]^{(\Omega f)_{\#}}
  & \pi_{n}(\Omega Y)
\ar[d]^{\partial_{\#}} \ar[r]^{j} &
  \pi_{n}(\Omega f)\ar[d]^{\Delta} \ar[r]^-
{p} & \cdots
    \\
\cdots\ar[r]^-{\widehat{p}} &
\pi_{n}(\map_*(X,X;1))
     \ar[r]^-{(f_*)_{\#}}  &
     \pi_{n}(\map_*(X,Y;f)) \ar[r]^-
{\widehat{j}} &
     \pi_{n}(f_*)\ar[r]^-{\widehat{p}} &
\cdots }
$$
Note that the maps in the $G$-sequence are
just the restrictions
of the maps in the long exact homotopy
sequence of the map $\Omega
f\colon \Omega X \to \Omega Y$. Thus
compositions of consecutive
maps in the $G$-sequence are trivial.
However, the sequence of
kernels of a commutative ladder of exact
sequences need not be
exact, and so the $G$-sequence is a chain
complex (of
$\Z$-modules).  The original description
given in \cite{L-W1,
L-W4} (see also \cite[Sec.1]{L-S})
represents the $G$-sequence as
an image sequence, in a way obviously
equivalent to the above.

Below, we construct a chain complex
associated to any DG Lie
algebra map that we will show corresponds to
the rationalized
$G$-sequence.  Observe that the relative
chain complex is a
functorial construction.  That is, given a
commutative square
$$\xymatrix{ V\ar[d]^{\alpha} \ar[r]^{\phi}
& W \ar[d]^{\beta} \\
V' \ar[r]^{\phi'} & W'}
$$
of DG vector space maps we obtain a DG
vector space map $(\beta,
\alpha)\colon \Rel(\phi) \to \Rel(\phi')$
given by $(\beta,
\alpha)(w, v) = (\beta(w), \alpha(v))$. This
leads to a
commutative diagram of short exact sequences
$$\xymatrix{0 \ar[r] & W_*
\ar[d]^{\beta}\ar[r]^-{J} &
\Rel_*(\phi)\ar[d]^{(\beta, \alpha)}
\ar[r]^-{P} &
V_{*-1} \ar[d]^{\alpha}\ar[r]& 0\\
0 \ar[r] & W'_* \ar[r]^-{J} & \Rel_*(\phi')
\ar[r]^-{P} & V'_{*-1}
\ar[r]& 0,}$$
and hence to a commutative ladder of long
exact homology
sequences.

Now suppose given a map $\psi \colon L \to
K$ of DG Lie algebras.
Apply the above to the commutative square
\begin{equation}\label{eq:commutative
square}
\xymatrix{ L \ar[d]^{\ad} \ar[r]^{\psi} & K
\ar[d]^{\ad_{\psi}} \\
\Der(L,L;1) \ar[r]^-{\psi_*} & \Der(L, K;
\psi)}
\end{equation}
of DG vector spaces.  Note that the relative
term $\Rel_*(\psi_*)$
is given by
$$\Rel_n(\psi_*) = \Der_n(L, K; \psi) \oplus
 \Der_{n-1}(L,L;1)$$
with differential $\delta(\theta_1,
\theta_2) = (\psi\circ\theta_2
- D(\theta_1), D(\theta_2))$. We obtain the
following commutative
ladder of long exact homology sequences:
\begin{equation}\label{eq:homology ladder of
psi}
\xymatrix{ \cdots \ar[r]^-{H(P)} & H_{n}(L)
\ar[d]^{H(\ad)}\ar[r]^{H(\psi)}  &  H_{n}(K)
\ar[d]^{H(\ad_{\psi})} \ar[r]^-{H(J)} &
H_n(\Rel(\psi))
\ar[d]^{H(\ad_{\psi}, \ad)} \ar[r]^-{H(P)}
&\cdots  \\
\cdots\ar[r]^-{H(\widehat{P})} &
H_n(\Der(L,L;1))
\ar[r]^-{H(\psi_*)}  & H_n(\Der(L,K;\psi))
\ar[r]^-{H(\widehat{J})} &
H_n(\Rel(\psi_*))\ar[r]^-
{H(\widehat{P})}&\cdots }
\end{equation}
To obtain unambiguous notation, we have
written $\widehat{J}\colon
\Der_n(L, K; \mathcal{L}_f) \to
\Rel_n(\psi_*)$ and
$\widehat{P}:\Rel_n(\psi_*) \to \Der_{n-
1}(L)$ for the usual
inclusion and projection maps in the lower
sequence.  We refer to
this ladder as the (horizontal) homology
ladder arising from the
diagram (\ref{eq:commutative square}).

\begin{theorem}\label{thm:equivalent
ladders}
Let $f\colon X \to Y$ be a map between
simply connected CW
complexes of finite type, with $X$ finite.
The rationalization of
the homotopy ladder (\ref{eq:homotopy
ladder}), down to the rung
$\partial_\#\otimes1 \colon \pi_2(\Omega
Y)\otimes\Q \to
\pi_{2}(\map_*(X,Y;f))\otimes\Q$, is
equivalent to the homology
ladder arising from
\begin{equation}\label{eq:Quillen ad
commutative square}
\xymatrix{ \mathcal{L}_X \ar[d]^{\ad}
\ar[r]^{\mathcal{L}_f} & \mathcal{L}_Y
\ar[d]^{\ad_{\mathcal{L}_f}} \\
\Der(\mathcal{L}_X,\mathcal{L}_X;1) \ar[r]^-
{(\mathcal{L}_f)_*} &
\Der(\mathcal{L}_X, \mathcal{L}_Y;
\mathcal{L}_f), }
\end{equation}
down to the rung $H(\ad_{\mathcal{L}_f})
\colon H_2(\mathcal{L}_Y)
\to
H_2(\Der(\mathcal{L}_X,\mathcal{L}_Y;\mathcal{L}_f))$.
\end{theorem}

\begin{proof}
Let $\Phi^Y\colon
\pi_n(\map_*(X,Y;f))\otimes\Q \to
H_n(\Der(\mathcal{L}_X,\mathcal{L}_Y;
\mathcal{L}_f))$ be the
isomorphism defined in the proof of
\thmref{thm:isos Phi and
Phi*}.  Let $\Phi^X$ be the isomorphism
obtained in the same way,
by specializing to the case in which $Y = X$
and $f = 1$.

Consider the (vertical) homotopy ladder
arising from
(\ref{eq:fibre connecting maps}), with exact
columns the long
exact sequences of the evaluation fibration
sequences $\Omega X
\to  \map_*(X,X;1) \to  \map(X,X;1)$ and
$\Omega Y \to
\map_*(X,Y;f) \to \map(X,Y;f)$.  From
\thmref{theorem:evaluation
sequence} applied to each column, we obtain
an equivalence between
this ladder and the \emph{vertical} homology
ladder arising from
(\ref{eq:Quillen ad commutative square}).
It follows that, with
the above notation, we have commutative
cubes
\begin{displaymath}
\xymatrix@C=-32pt{ & H_{n}(\mathcal{L}_{X})
 \ar[rr]^-{H(\mathcal{L}_{f})}
\ar'[d]^-{H(\ad)}[dd] & &
H_{n}(\mathcal{L}_{Y})
\ar[dd]^-{H(\ad_{\mathcal{L}_f})} \\
\pi_{n}(\Omega X)\otimes\Q
\ar[ru]^{\alpha_X}_{\cong}\ar[rr]_(0.63){
(\Omega f)_\#\otimes1}
 \ar[dd]_{\partial_\#}& &
\pi_{n}\big(\Omega Y\big)\otimes\Q
\ar[ru]^{\alpha_Y}_{\cong}
\ar[dd]^(0.3){\partial_\#}\\
 &
H_{n}\big(\Der(\mathcal{L}_{X},\mathcal{L}_{X};1)\big)
\ar'[r]_(0.9){H((\mathcal{L}_{f})_*)}[rr] &
&
H_{n}\big(\Der(\mathcal{L}_{X},\mathcal{L}_{Y};\mathcal{L}_{f})\big)  \\
\pi_n\big(\map_*(X, X; 1)\big)\otimes\Q
\ar[ru]^{\Phi^X}_{\cong}
\ar[rr]_{(f_*)_\#\otimes1} & &
\pi_n\big(\map_*(X, Y;f)\big)\otimes\Q\ar[ru]^{\Phi^Y}_{\cong} }
\end{displaymath}
for $n\geq2$.  Here, $\alpha_X$ and
$\alpha_Y$ denote the
isomorphisms that are obtained from
\lemref{lemma:alpha prime}.
Strictly speaking, we cannot conclude an
equivalence of ladders
from \thmref{theorem:evaluation sequence},
due to the non-natural
choices made in \lemref{lemma:alpha prime}.
However, an easy
extension of \lemref{lemma:alpha prime} to
the setting of an
equivalence of ladders overcomes this
difficulty---see
\cite[Lem.3.9]{L-S} for details.

Switching now to horizontal ladders, we may
use top and bottom
faces of the displayed cubes, together with
the extension of
\lemref{lemma:alpha prime} referred to
already, to obtain the
desired equivalence of ladders.
\end{proof}

\begin{remark}
Notice that \thmref{thm:equivalent ladders}
contains a description
of the long exact rational homotopy
sequences of a general map
$f$, and of the induced map $f_*$.  See
\cite[Th.3.3, Th.3.5]{L-S}
for the corresponding descriptions of these
sequences in terms of
Sullivan minimal models.
\end{remark}

In \defref{definition:evaluation group}, we
have defined
evaluation subgroups of a map of DG Lie
algebras, and of a DG Lie
algebra.  Our cast of characters that appear
in the $G$-sequence
is completed by the following:

\begin{definition}\label{definition:relative
evaluation group}
Let $\psi\colon L \to K$ be a DG Lie algebra
map. The \emph{$n$th
relative evaluation subgroup} $G_n^{rel}(K,
L; \psi)$ of $\psi$ is
the subgroup
$$G^{rel}_n(K, L; \psi) = \ker \{
H(\ad_{\psi}, \ad)\colon H_{n-1}(\Rel(\psi))
\to H_{n-1}(\Rel(\psi_*)) \}$$
of $H_n(\Rel(\psi))$.   The \emph{$G$-sequence of $\psi$} is the
sequence of kernels from the commutative
ladder (\ref{eq:homology
ladder of psi}).  That is, the sequence
$$
\xymatrix{  \cdots \ar[r]^-{H(P)} & G_{n}(L)
\ar[r]^-{H(\psi)} &
G_n(K,L;\psi) \ar[r]^-{H(J)} & G_n^{rel}(K,
L; \psi)
\ar[r]^-{H(P)} &\cdots.}
$$
\end{definition}

\begin{corollary}\label{cor:G-sequence}
Let $f\colon X \to Y$ be a map between
simply connected CW
complexes of finite type, with $X$ finite.
Then the
rationalization of the $G$-sequence of $f$
$$\xymatrix{\cdots \ar[r] &  G_{n}(X)
\ar[r]^-{\Omega f_\#\otimes\Q}  & G_{n}(Y,
X;f) \otimes \Q
\ar[r]^{j} & G_{n}^{rel}(Y, X;f) \otimes \Q
\ar[r] &  \cdots}$$
down to the term $G_{3}(Y, X;f)\otimes\Q$ is
equivalent to the
$G$-sequence of the Quillen model
$\mathcal{L}_f\colon
\mathcal{L}_X \to \mathcal{L}_Y$ of $f$,
$$\xymatrix{\cdots \ar[r] &
G_{n}(\mathcal{L}_X)
\ar[r]^-{H(\mathcal{L}_f)} &
G_n(\mathcal{L}_Y,
\mathcal{L}_X;\mathcal{L}_f) \ar[r]^{H(J)} &
G_n^{rel}(\mathcal{L}_Y, \mathcal{L}_X;
\mathcal{L}_f)  \ar[r] &
\cdots }$$
down to the term  $G_3(\mathcal{L}_Y,
\mathcal{L}_X;
\mathcal{L}_f)$.
\end{corollary}

\begin{proof}
This follows immediately from
\thmref{thm:equivalent ladders},
since equivalent ladders have equivalent
kernel (and image)
sequences.
\end{proof}

In particular, we can add the following to
\thmref{theorem:evaluation group}

\begin{corollary}\label{cor:G rel terms}
Let $f\colon X \to Y$ be a map between
simply connected CW
complexes of finite type with $X$ finite.
Then the rationalized
relative evaluation subgroup $G^{rel}_n(Y,
X;f)$ is isomorphic to
the relative Gottlieb group of the Quillen
model
$\mathcal{L}_f\colon \mathcal{L}_X \to
\mathcal{L}_Y$ of $f$ as in
\defref{definition:relative evaluation
group}. That is, for $n \geq 4$
$$G^{rel}_n(Y,X;f) \otimes \Q \cong   \ker\{
H(\ad_{\mathcal{L}_f}, \ad)\colon H_{n-
1}(\Rel(\mathcal{L}_f))  \to H_{n-1}(\Rel(
\mathcal{L}_f^*))
\}.$$
\end{corollary}

\begin{remark}
With a little more work, the equivalence of
\corref{cor:G-sequence}, and hence
\corref{cor:G rel terms}, may
be extended to the $G^{rel}_3(Y, X;f)$ term.
Under our
hypotheses, $G_2(X)\otimes\Q = 0$, and there
seems little to be
gained by trying to extend \corref{cor:G-sequence} beyond this
point (cf.~the remarks on the low-end terms
in
\cite[Rem.3.11]{L-S}).
\end{remark}

As an application of the above, we consider
the question of
exactness of the $G$-sequence for cellular
extensions. The first
example of inexactness for a cellular
inclusion of finite
complexes was given in \cite{L-W2}. Observe
that rational examples
of non-exactness delocalize, since exactness
in the ordinary
setting implies exactness after
rationalization.

We focus on the rationalized $G$-sequence at
the $G_n(X)$ term.
Following Lee and Woo \cite{L-W4}, define
the $\omega$-homology of
a map $f\colon X \to Y$ at this term by
setting
$$H_n^{a\omega}(Y, X; f) = \frac{ \ker \{
f_\#\colon G_n(X) \to
G_n(Y,X;f) \} }{\im \{ p\colon
G_{n+1}^{rel}(Y,X;f) \to G_n(X) \}
}.$$
For a single cell-attachment, the following
is a complete result
for the lowest degree in which the
(rationalized) $G$-sequence can
be non-exact at the $G_n(X)$-term.

\begin{theorem}\label{thm:G-Sequence one cell}
Let $X$ be a finite CW
complex of dimension $n$, and $Y = X
\cup_\alpha e^{m+1}$ for some
$\alpha \in \pi_m(X)$.  Suppose the
following three conditions
hold:
\begin{enumerate}
\item  $\alpha_\Q \not = 0$;
\item $\alpha_\Q \in G_m(X_\Q)$;
\item $Y$ is not rationally equivalent to a
point.
\end{enumerate}
Then the $G$-sequence of the inclusion $i
\colon X \to Y$ is
non-exact at the $G_m(X)$ term, and $H_m^{a
\omega}(Y,X;i) \otimes
\Q = \Q$. Conversely, if any of (1)--(3) do
not hold, then $H_m^{a
\omega}(Y,X;i) \otimes \Q = 0$, that is,
$$\xymatrix{G_{m+1}^{rel}(Y,X;i)\otimes\Q
\ar[r]^-{p\otimes\Q}& G_m(X)\otimes\Q
\ar[r]^-{i_\#\otimes\Q}
& G_{m}(Y,X;i)\otimes\Q}$$
is exact.
\end{theorem}

\begin{proof}
Clearly, the kernel of $i_\#\otimes\Q \colon
\pi_m(X_\Q) \to
\pi_m(Y_\Q)$ is the subspace
$\langle\alpha_\Q\rangle$ of
$\pi_m(X_\Q)$. If $\alpha_\Q =0$, then
$i_\#\otimes\Q :\pi_m(X_\Q)
\to \pi_m(Y_\Q)$ is injective, and if
$\alpha_\Q \not\in
G_m(X_\Q)$, then the restriction of
$i_\#\otimes\Q$ to $G_n(X_\Q)$
of $i_\#$ is injective.  In either case,
$H_m^{a \omega}(Y,X;i)
\otimes \Q = 0$ by definition.

So suppose that (1) and (2) hold.  We show
that
$G_{m+1}^{rel}(Y,X;i)\otimes\Q = 0$ if (3)
also holds. For this,
we use \corref{cor:G rel terms} and show
that
$H(\ad_{\mathcal{L}_i}, \ad)\colon
H_{m}(\Rel(\mathcal{L}_i))  \to
H_{m}(\Rel( (\mathcal{L}_i)_*))$ is
injective.

We may write the Quillen model of $i\colon X
\to Y$ as an
inclusion $\mathcal{L}_i\colon \mathcal{L}_X
\to \mathcal{L}_Y$,
with $\mathcal{L}_X = \L(W)$, $\mathcal{L}_Y
= \L(W)\sqcup\L(y)$,
and differentials $d_Y | \L(W) = d_X$ and
$d_Y(y) = \chi \in
\L(W)$ a cycle of degree $(m-1)$ in
$\mathcal{L}_X$ whose homology
class represents $\alpha_\Q$ via the
isomorphism
$H_{m-1}(\mathcal{L}_X) \cong
\pi_m(X)\otimes\Q$
(cf.~\cite[Sec.24(d)]{F-H-T}).

Any cycle $\zeta \in \Rel_m(\mathcal{L}_i) =
(\mathcal{L_Y})_m
\oplus (\mathcal{L}_X)_{m-1}$ may be written
$\zeta = (\lambda y +
\xi, d_Y(\lambda y + \xi))$ for
$\lambda\in\Q$ and $\xi \in
(\mathcal{L}_X)_m$.  Suppose that
$H(\ad_{\mathcal{L}_i},
\ad)(\langle\zeta\rangle) = 0$.  Then
$(\ad_{\mathcal{L}_i},
\ad)(\zeta) = \delta(\varphi, \theta)$ for
some $(\varphi, \theta)
\in \Rel_{m+1}((\mathcal{L}_i)_*)$.  In
particular, we have
$\ad_{\mathcal{L}_i}(\lambda y + \xi) =
(\mathcal{L}_i)_*(\theta)
- D\varphi \in
\Der_m(\mathcal{L}_X,\mathcal{L}_Y;\mathcal{
L}_i)$.
We now use (3), if necessary, to choose an
indecomposable $w \in
W$ such that $\chi$ and $w$ are linearly
independent in
$\mathcal{L}_X$. (A choice is only necessary
when $d_Y(y)$ is
indecomposable.  If $\chi$ is decomposable,
then any
indecomposable of $W$ will do.)  On this
indecomposable, we
evaluate as follows:
\begin{equation}\label{eq:left-side}
\ad_{\mathcal{L}_i}(\lambda y + \xi)(w) =
\lambda [y, w] + [\xi,
w]
\end{equation}
and
\begin{equation}\label{eq:right-side}
((\mathcal{L}_i)_*(\theta) - D\varphi)(w)  =
\mathcal{L}_i\circ\theta(w) - d_Y\varphi(w)
+ (-1)^{m+1}\varphi
d_X(w)
\end{equation}
We claim that all terms of (\ref{eq:right-side}) are independent
of $[y,w]$.  First,
$\mathcal{L}_i\circ\theta(w) \in \L(W)$.
Next,
$\varphi(w) \in
(\L(W)\sqcup\L(y))_{m+1+|w|}$ may contain
terms in
$\L(W)$, quadratic terms $[y, w']$ for $w'
\in W$ or $[y,y]$ (if
$m$ is odd), or terms involving $y$ of
bracket-length at least
three. Now $d_Y([y, w']) = [\chi, w'] \pm
[y, d_Y(w')]$, and
$d_Y([y,y]) = 2 [\chi, y]$ (if $m$ is odd).
By choice, $\chi$ and
$w$ are linearly independent, and also $d_X$
is decomposable.
Hence, all terms of $d_Y\varphi(w)$ are
linearly independent of
$[y,w]$.  Finally, $d_X(w)$ is decomposable
and in
$\L(W)_{|w|-1}$.  Hence, $\varphi d_X(w)$ is
in the ideal of
$\L(W)\sqcup\L(y)$ generated by elements of
$W$ of degree $\leq
|w|-2$.  This proves the claim, and shows
that if
(\ref{eq:left-side}) and (\ref{eq:right-side}) agree, which must
be the case if $H(\ad_{\mathcal{L}_i},
\ad)(\langle\zeta\rangle) =
0$, then we have $\lambda = 0$.  We have
shown that the kernel of
$H(\ad_{\mathcal{L}_i}, \ad)\colon
H_{m}(\Rel(\mathcal{L}_i))  \to
H_{m}(\Rel( (\mathcal{L}_i)_*))$ consists of
classes represented
by cycles of the form $\zeta = (\xi,
d_Y(\xi))$ for $\xi \in
(\mathcal{L}_X)_m$.  Since $d_Y(\xi) =
d_X(\xi)$, we may write
$\delta(0, \xi) = (\xi, d_X(\xi) = \zeta$.
That is,
$H(\ad_{\mathcal{L}_i}, \ad)\colon
H_{m}(\Rel(\mathcal{L}_i))  \to
H_{m}(\Rel( (\mathcal{L}_i)_*))$ is
injective and hence by
\corref{cor:G rel terms},
$G_{m+1}^{rel}(Y,X;i)\otimes\Q = 0$.
This establishes the first assertion of the
theorem.

The second assertion will be established
when we handle the case
in which (1) and (2) hold, but (3) does not,
that is, in which $Y
\simeq_\Q *$. In this case, we must have $X
\simeq_\Q S^n$, $m=n$,
and $\alpha \in \pi_n(S^n)$ non-trivial. The
Quillen model of $Y$
is then $\L(w, y)$ with $|w| = n-1$, $|y| =
n$, and $d_Y(y) = w$.
(Notice that this is precisely the case in
which we cannot choose
a second indecomposable of $W$ independent
of $\chi$.)  In
addition, since we are assuming that (2)
holds, we must have $m$
odd, and hence $|w|$ is even and $\L(w)$ is
an abelian Lie
algebra. It is easy to see that $(y,w) \in
\Rel_m(\mathcal{L}_i)$
is a cycle that cannot be a $\delta$-boundary. Define an
$\mathcal{L}_i$-derivation $\varphi \in
\Der_{m+1}(\mathcal{L}_X,\mathcal{L}_Y;\mathcal{L}_i)$ by setting
$\varphi(w) = -\frac{1}{2} [y,y]$.  A
straightforward check shows
that $(\ad_{\mathcal{L}_i}, \ad)(y,w) =
(\ad_{\mathcal{L}_i}(y),
0) = \delta(\phi, 0)$ in $\Rel_m(
(\mathcal{L}_i)_*)$.  That is,
$\langle y, w \rangle$ represents an element
in
$G^{rel}_m(\mathcal{L}_Y,\mathcal{L}_X;\mathcal{L}_i)$.
Furthermore, $H(P)(\langle y, w \rangle) =
\langle w \rangle$, and
so $H(P) \colon
G^{rel}_m(\mathcal{L}_Y,\mathcal{L}_X;\mathcal{L}_i) \to
G_m(\mathcal{L}_X)$ is onto.  The result
follows.
\end{proof}

We conclude with an example of vanishing
rational
$\omega$-homology.  In the following result,
we use the ideas
discussed before \thmref{thm:coformal G=P},
concerning the notion
of a coformal map.

\begin{theorem}\label{thm:G-sequence exact
coformal map}
Let $f\colon X \to Y$ be a coformal map
between CW complexes of
finite type, with $X$ finite.  Then
$H_n^{a\omega}(Y, X; f)
\otimes \Q = 0$, that is,
$$\xymatrix{G_{n+1}^{rel}(Y,X;f)\otimes\Q
\ar[r]^-{p\otimes\Q}& G_n(X)\otimes\Q
\ar[r]^-{f_\#\otimes\Q}
& G_{n}(Y,X;f)\otimes\Q}$$
is exact, for each $n \geq 3$.
\end{theorem}

\begin{proof} We will use \corref{cor:G-sequence} and show that
$\ker \{ H(\mathcal{L}_f) \} \subseteq \im
\{ H(P)\}$. From
\defref{definition:coformal map}, we assume
that both
$\mathcal{L}_X$ and $\mathcal{L}_Y$ admit
upper (second) gradings
with the properties described in
\defref{definition:coformal Lie
algebra}, and that $\mathcal{L}_f$ preserves
upper degrees.  Let
$\alpha = \langle \xi \rangle \in
G_n(\mathcal{L}_X)$ satisfy
$H(\mathcal{L}_f)(\alpha) = 0$.  We assume
$\xi$ is of upper
degree zero in the bigraded model for
$\mathcal{L}_X$.
Furthermore, $\mathcal{L}_f(\xi) = d_Y(y)$
for some $y \in
(\mathcal{L}_Y )_{n+1}$ that we may assume
is of upper degree $1$.
Since $\alpha$ is Gottlieb,  $\ad(\xi) =
D(\psi)$ for some
derivation $\psi \in
\Der_{n+1}(\mathcal{L}_X, \mathcal{L}_X;1)$
and using the coformality of $\mathcal{L}_X$
again, we may assume
$\psi$ increases upper degree homogeneously
by $1$. The pair $(y,
\xi) \in \Rel_{n+1}(\mathcal{L}_f)$ is a
$\delta$-cycle that
satisfies $P(y, \xi) = \xi$. We now show
that $(y, \xi)$
represents an element in
$G_{n+1}^{rel}(\mathcal{L}_Y,\mathcal{L}_X;\mathcal{L}_f)$, that
is, we show the pair
$(\ad_{\mathcal{L}_f}(y), \ad(\xi))$ bounds
in $\Rel_{n+1}((\mathcal{L}_f)_*)$.   Set
$\Theta =
(\mathcal{L}_f)_*(\psi) -
\ad_{\mathcal{L}_f}(y)$, a derivation in
$\Der_{n+1}(\mathcal{L}_Y, \mathcal{L}_X;
\mathcal{L}_f)$. It is
direct to check that $D(\Theta) = d_Y \circ
\Theta
-(-1)^{n+1}\Theta \circ d_X = 0$. Moreover,
$\Theta$ increases
upper degree homogeneously by $1$.  Now
adapt the proof of
\thmref{thm:coformal G=P} to the current
situation, by replacing
the derivation $\ad_{\mathcal{L}_f}(\xi)$ in
that proof by $\psi$.
The inductive argument used there now
results in a derivation
$\theta\in \Der_{n+2}(\mathcal{L}_X,
\mathcal{L}_Y;
\mathcal{L}_f)$, constructed in the same way
only increasing upper
degree by $2$, that satisfies $\Theta =
D(\theta)$. Then the pair
$(\theta, \psi) \in
\Rel_{n+2}((\mathcal{L}_f)_*)$ satisfies
$\delta(\theta, \psi) =
(\ad_{\mathcal{L}_f}(y), \ad(\xi))$.
\end{proof}

\begin{appendix}

\section{Some DG Lie Algebra Homotopy
Theory}\label{section:DG Homotopy Theory}

In this appendix, we present some DG Lie
algebra homotopy theory.
Our main references for this material are
\cite[Ch.II.5]{Tan83}
and part IV of \cite{F-H-T}.  For
completeness and convenience, we
recall the basic notions here.  The main
point of the appendix is
to provide details for some facts used in a
crucial way to
establish our main results.  Since this is a
technical appendix,
we assume a greater degree of familiarity
with techniques from
rational homotopy theory than in the main
body of the paper.

The algebraic notion of homotopy that we use
here is ``left
homotopy" of DG Lie algebra maps, defined in
terms of a suitable
cylinder object for a free DG Lie algebra
$\L(V;d)$.  This is
another free DG Lie algebra denoted
$\L(V)_I$, together with
inclusions $\lambda_0, \lambda_1 \colon
\L(V) \to \L(V)_I$ and a
projection $p\colon \L(V)_I \to \L(V)$ that
together satisfy
$p\circ\lambda_i = 1$ for $i = 0,1$.  As a
DG Lie algebra, we have
$\L(V)_I = \L(V, sV, \widehat{V}; D)$, with
$sV$ the suspension of
$V$, and $\widehat{V}$ an isomorphic copy of
$V$. The differential
$D$ in $\L(V)_I$ extends $d$, so that
$\L(V)$ is a sub-DG Lie
algebra, and is defined on the other
generators as $D(sv) =
\widehat{v}$ and $D(\widehat{v}) = 0$, for
$v \in V$. The
inclusion $\lambda_0 \colon \L(V) \to
\L(V)_I$ is the obvious
inclusion of the sub-DG Lie algebra $\L(V)$.
The inclusion
$\lambda_1 \colon \L(V) \to \L(V)_I$, on the
other hand, is more
involved.  Define a derivation $\sigma \in
\Der_1(\L(V)_I,
\L(V)_I; 1)$ on generators by $\sigma(v) =
sv$, $\sigma(sv) = 0$,
and $\sigma(\widehat{v}) = 0$, then define
$\theta \in
\Der_0(\L(V)_I, \L(V)_I; 1)$ as $\theta =
[D, \sigma] = D\circ
\sigma + \sigma\circ D$.  Observe that
$\theta$ is a cycle, and is
locally nilpotent.  Therefore,
exponentiating $\theta$ gives a
(DG) automorphism $\exp(\theta) \colon
\L(V)_I \to \L(V)_I$.
Finally, define $\lambda_1 =
\exp(\theta)\circ\lambda_0$.  The
projection $p$ is defined in the obvious way
as $p(v) = v$, $p(sv)
= 0$, and $p(\widehat{v}) = 0$. Given a pair
of DG Lie algebra
maps $\phi, \psi\colon \L(V) \to L$, we say
\emph{$\phi$ is
homotopic to $\psi$} if there exists a map
$\mathcal{H}\colon
\L(V)_I \to L$ such that
$\mathcal{H}\circ\lambda_0 = \phi$ and
$\mathcal{H}\circ\lambda_1 = \psi$.  In this
case, we say
$\mathcal{H}$ is a (DG) homotopy from $\phi$
to $\psi$.

In addition to the notation established
above, which we take as
fixed for this appendix, we will use the
following conventions for
DG Lie algebra maps: We will generally
suppress differentials from
our notation.  Recall the model for
$S^n\times X$ from
\secref{section:Quillen model}.  We will use
$J$ to denote either
inclusion $\L(W) \to \L(W, v, W')$ or $\L(W,
v) \to \L(W, v, W')$,
and $J_I$ to denote the corresponding
inclusions of cylinder
objects.  Namely, in the first case, $J_I(w)
= w$, $J_I(sw) = sw$,
and $J_I(\widehat{w}) = \widehat{w}$, for $w
\in W$. We also fix
some notation for maps of spaces that we use
in this appendix. Let
$i_1 \colon S^n \to S^n \times X$ and $i_2
\colon X \to S^n \times
X$ denote the inclusions, and $\pi_1$ the
projection onto the
first factor of a product of spaces.

The basic correspondence between the notions
of homotopy in the
topological and algebraic settings is as
follows: Maps $f, g\colon
X \to Y$ of rational spaces are homotopic if
and only if their
Quillen models $\mathcal{L}_{f},
\mathcal{L}_{g}\colon
\mathcal{L}_{X} \to \mathcal{L}_{Y}$ are
homotopic in the sense
just defined. For the proof of
\thmref{thm:isos Phi and Phi*},
however, we need finer detail than this
basic correspondence.
Consider the argument to show $\Phi'$ is
well-defined, for
example.  Suppose $a, b \colon S^n \to
\map_*(X,Y;f)$ are
homotopic representatives of the same
homotopy element, and that
their adjoints are $A, B \colon S^n \times X
\to Y$.  Then $A$ and
$B$ are homotopic relative to $S^n \vee X$.
Indeed, if $H$ denotes
the homotopy from $A$ to $B$ that is adjoint
to the homotopy from
$a$ to $b$, then the following diagram
commutes:
\begin{equation}\label{eq:based homotopy}
\xymatrix{(S^n\vee X)\times I
\ar[d]_{(i_1|i_2)\times 1}
\ar[r]^-{\pi_1} &
S^n \vee X \ar[d]^{({*}\mid f)} \\
S^n\times X \times I \ar[r]_-{H} & Y}
\end{equation}
For our argument that $\Phi'$ is well-
defined, it is crucial that
we may assume the DG Lie algebra homotopy
that corresponds to $H$
is restricted in a certain way.
Specifically, if $\mathcal{H}$
denotes the DG Lie algebra homotopy that
corresponds to $H$, then
we require that the following diagram
commute:
\begin{equation}\label{eq:based DG Lie homotopy}
\xymatrix{\L(W, v)_I  \ar[d]_{J_I} \ar[r]^-
{p} &
\L(W, v) \ar[d]^{(\mathcal{L}_f\mid 0)} \\
\L(W, v, W')_I \ar[r]_-{\mathcal{H}} &
\mathcal{L}_Y}
\end{equation}
We also need a similar fact for homotopy
elements of the unbased
mapping space to establish that $\Psi'$ is
well-defined.  Whilst
this translation is intuitively plausible,
there are some
technical details to be checked.

Our starting point for this is the
corresponding result in the
Sullivan model setting, which we have proved
in our earlier paper.
We assume familiarity with the usual
notation in that setting,
that is, rational homotopy theory from the
DG algebra point of
view (see \cite{F-H-T}). The algebraic
notion of homotopy that we
used in \cite{L-S} is ``right homotopy"
defined in terms of a
suitable path object for a free DG algebra.
For $\Lambda(V)$, this
consists of maps
\begin{equation}\label{eq:path in DGA}
\xymatrix@C=30pt{\Lambda(V) \ar[r]^-{j} &
\Lambda(V)\otimes\Lambda(t, dt) \ar[r]^-
{p_0, p_1} & \Lambda(V)}
\end{equation}
that satisfy $p_i\circ j = 1$ for $i = 0,1$.
Here, $j$ is the
inclusion $j(\chi) = \chi\otimes1$.  The
projections are defined
by $p_i(t) = i$, $p_i(dt) = 0$, and
$p_i(\chi) = \chi$ for $\chi
\in \Lambda(V)$. Given a pair of DG algebra
maps $\phi, \psi\colon
\Lambda(W) \to \Lambda(V)$, we say
\emph{$\phi$ is homotopic to
$\psi$} if there exists a map
$\mathcal{G}\colon \Lambda(W) \to
\Lambda(V)\otimes\Lambda(t, dt)$ such that
$p_0\circ\mathcal{G} =
\phi$ and $p_1\circ\mathcal{G} = \psi$.  In
this case, we say
$\mathcal{G}$ is a (DG) homotopy from $\phi$
to $\psi$.

We fix more notation: Let $\mathcal{M}_X$
denote the Sullivan
minimal model of a space $X$.  Recall that
the minimal model of a
product of spaces is the tensor product of
their minimal models,
so that
$\mathcal{M}_{S^n}\otimes\mathcal{M}_X$ is a
Sullivan
model for $S^n \times X$.  Let
$\mathcal{M}_f \colon \mathcal{M}_Y
\to \mathcal{M}_X$ denote a Sullivan minimal
model of a map
$f\colon X \to Y$.  Let $q_1 \colon
\mathcal{M}_{S^n}\otimes\mathcal{M}_X \to
\mathcal{M}_{S^n}$ and
$q_2 \colon
\mathcal{M}_{S^n}\otimes\mathcal{M}_X \to
\mathcal{M}_{X}$ denote the projections.

\begin{lemma}\label{lem:restricted DG Alg
homotopy}
Let $f \colon X \to Y$ be a map with a fixed
choice of Sullivan
model $\mathcal{M}_f \colon \mathcal{M}_Y
\to \mathcal{M}_X$. Let
$A, B \colon S^n\times X \to Y$ be maps that
restrict to $A\circ
i_2 = B\circ i_2 = f \colon X \to Y$.

\begin{enumerate}
\item[(i)] Suppose $A \sim B$ via a homotopy
$H$ relative to $X$, that is,
suppose the diagram
\begin{equation}\label{eq:unbased homotopy}
\xymatrix{X\times I \ar[d]_{i_2\times 1}
\ar[r]^-{\pi_1} &
 X \ar[d]^{f} \\
S^n\times X \times I \ar[r]_-{H} & Y}
\end{equation}
commutes.  Then there exists a homotopy
$\mathcal{G} \colon
\mathcal{M}_Y \to
\mathcal{M}_{S^n}\otimes\mathcal{M}_X\otimes
\Lambda(t, dt)$ from
$\mathcal{M}_A$ to $\mathcal{M}_B$ that is
``relative to
$\mathcal{M}_X$," in the sense that the
diagram
\begin{equation}\label{eq:unbased Sullivan homotopy}
\xymatrix{\mathcal{M}_Y \ar[r]^-
{\mathcal{G}}
\ar[d]_{\mathcal{M}_f}&
\mathcal{M}_{S^n}\otimes\mathcal{M}_X\otimes
\Lambda(t, dt)
\ar[d]^{q_2\otimes1}\\
\mathcal{M}_X \ar[r]_-{j}&
\mathcal{M}_X\otimes\Lambda(t, dt) }
\end{equation}
(strictly) commutes.
\item[(ii)] Suppose further that $A$ and $B$
restrict to $A\circ
i_1 = B\circ i_1 = * \colon S^n \to Y$, and
the homotopy $H$ is
also relative to $S^n$, so that
(\ref{eq:based homotopy})
commutes.  Then there exists a homotopy
$\mathcal{G} \colon
\mathcal{M}_Y \to
\mathcal{M}_{S^n}\otimes\mathcal{M}_X\otimes
\Lambda(t, dt)$ from
$\mathcal{M}_A$ to $\mathcal{M}_B$ that is
``relative to
$\mathcal{M}_{S^n\vee X}$," in the sense
that the diagram
$$\xymatrix{\mathcal{M}_Y \ar[r]^-
{\mathcal{G}}
\ar[d]_{(\varepsilon,\mathcal{M}_f)}&
\mathcal{M}_{S^n}\otimes\mathcal{M}_X\otimes
\Lambda(t, dt)
\ar[d]^{(q_1\otimes1, q_2\otimes1)}\\
\mathcal{M}_{S^n}\oplus\mathcal{M}_X
\ar[r]_-{j\oplus j}&
\mathcal{M}_{S^n}\otimes\Lambda(t,
dt)\oplus\mathcal{M}_X\otimes\Lambda(t, dt)
}
$$
(strictly) commutes ($\varepsilon \colon
\mathcal{M}_{Y} \to
\mathcal{M}_{S^n}$ denotes the map that is
zero in positive
degrees).
\end{enumerate}
\end{lemma}

\begin{proof}
The first point is proved in
\cite[Lem.A.2]{L-S}.  The argument
given there is readily adapted to prove the
second point.
\end{proof}

We will carefully translate this result into
the Quillen model
setting via the so-called Quillen functor
$\QL$.  This is a
functor from the category of differential
graded algebras to the
category of differential graded Lie algebras
(see
\cite[Sec.22(e)]{F-H-T} or
\cite[I.1.(7)]{Tan83} for details).
Since we only use general properties of this
functor, we do not
recall its definition here. We do recall
that it preserves
quasi-isomorphisms and, as a consequence of
its definition, takes
an injective DG algebra map to a surjective
DG Lie algebra map.

Sullivan and Quillen models of a map are
related via the Quillen
functor.  Suppose $f\colon X \to Y$ has
Sullivan minimal model
$\mathcal{M}_f \colon \mathcal{M}_Y \to
\mathcal{M}_X$.  Applying
$\QL$ gives $\QL(\mathcal{M}_f) \colon
\QL(\mathcal{M}_X) \to
\QL(\mathcal{M}_Y)$.  Now suppose $\rho_X
\colon \mathcal{L}_X \to
\QL(\mathcal{M}_X)$ and $\rho_Y \colon
\mathcal{L}_Y \to
\QL(\mathcal{M}_Y)$ are given minimal
models. Then any map
$\mathcal{L}_f \colon \mathcal{L}_X \to
\mathcal{L}_Y$ such that
$\rho_Y\circ\mathcal{L}_f$ and
$\QL(\mathcal{M}_f)\circ\rho_X$ are
homotopic is a Quillen minimal model for
$f$. Recall that we
obtain such maps as follows: There is a
standard way of converting
the quasi-isomorphism  $\rho_Y$ into a
surjective
quasi-isomorphism.  Namely, let
$E(\QL(\mathcal{M}_Y))$ denote the
acyclic, free DG Lie algebra $\L(V, s^{-
1}V)$, with $V$ a vector
space isomorphic to $\QL(\mathcal{M}_Y)$ and
$d(v) = s^{-1}v$ for
$v \in V$.  Let $\widehat{\rho_Y} \colon
\mathcal{L}_Y\sqcup
E(\QL(\mathcal{M}_Y)) \to
\QL(\mathcal{M}_Y)$ be the map that
extends $\rho_Y$ on $\mathcal{L}_Y$, and
maps $\widehat{\rho_Y}(v)
= v$, $\widehat{\rho_Y}(s^{-1}v) = dv$ for
$v \in V$.  Then we may
lift $\QL(\mathcal{M}_f)\circ \rho_X$
through the surjective
quasi-isomorphism $\widehat{\rho_Y}$,
$$\xymatrix{\L(W) \ar[d]_{\rho_X}^{\simeq}
\ar@{.>}[r]^-{\phi_f} &
\mathcal{L}_Y\sqcup E(\QL(\mathcal{M}_Y))
\ar@{->>}[d]^{\widehat{\rho_Y}}_{\simeq}
\ar[r]^-{\pi}
& \mathcal{L}_Y \ar[ld]^{\rho_Y} \\
\QL(\mathcal{M}_X) \ar[r]_-
{\QL(\mathcal{M}_f)} &
 \QL(\mathcal{M}_Y) }
$$
to obtain a map $\phi_f$ that satisfies
$\widehat{\rho_Y}\circ\phi_f =
\QL(\mathcal{M}_f)\circ\rho_X$ (see
\cite[II.5.(13)]{Tan83} or
\cite[Prop.22.11]{F-H-T} for the
standard lifting lemma). A Quillen minimal
model for $f$ is then
obtained by composing with the projection
$\pi$ to give
$\mathcal{L}_f = \pi\circ\phi_f$.

Homotopies are also related via the Quillen
functor.  However,
since it is contravariant, the Quillen
functor translates a right
homotopy of DG algebra maps into a left
homotopy of DG Lie algebra
maps, in the following way. By applying
$\QL$ to (\ref{eq:path in DGA}), we obtain maps
$$\xymatrix@C=40pt{\QL(\mathcal{M}_X)
\ar[r]^-{\QL(p_0), \QL(p_1)}  &
\QL\big(\mathcal{M}_X\otimes\Lambda(t,
dt)\big) \ar[r]^-{\QL(j)} &
\QL(\mathcal{M}_X)}$$
that satisfy $\QL(j)\circ \QL(p_i) = \QL(1)
= 1$ for $i= 0,1$.
Since $j$ is an injective quasi-isomorphism,
$\QL(j)$ is a
surjective quasi-isomorphism.   Since
$\QL(j)\circ\QL(p_0)\circ\rho_X = \rho_X =
\QL(j)\circ\QL(p_1)\circ\rho_X \colon
\mathcal{L}_X \to
\QL(\mathcal{M}_X)$, the linear difference
$\QL(p_1)\circ\rho_X -
\QL(p_0)\circ\rho_X \colon \mathcal{L}_X \to
\QL(\mathcal{M}_X\otimes\Lambda(t, dt))$ has
image contained in
the DG ideal $\ker(\QL(j))$ of
$\QL(\mathcal{M}_X\otimes\Lambda(t,
dt))$.  Now suppose $\mathcal{L}_X = \L(W)$.
Since $\QL(j)$ is a
surjective quasi-isomorphism, $\ker(\QL(j))$
is an acyclic DG
ideal.  Hence, by a standard argument
(\cite[Prop.II.5(4)]{Tan83}), we may
construct a left-homotopy $G
\colon \L(W)_I \to
\QL(\mathcal{M}_X\otimes\Lambda(t, dt))$
from
$\QL(p_0)\circ\rho_X$ to
$\QL(p_1)\circ\rho_X$ \emph{that in
addition satisfies} $G\big(\L(sW,
\widehat{W})\big) \subseteq
\ker(\QL(j))$ (this last point is key). That
is, we have the
following commutative diagram,
$$\xymatrix{\L(W)
\ar[rrd]^{\QL(p_0)\circ\rho_X}
\ar[d]_{\lambda_0} \\
\L(W)_I \ar[rr]^(0.3){G} & &
\QL(\mathcal{M}_X\otimes\Lambda(t,
dt))  \\
\L(W)\ar[rru]_{\QL(p_1)\circ\rho_X}
\ar[u]^{\lambda_1} }
$$
for which $G\big(\L(sW, \widehat{W})\big)
\subseteq \ker(\QL(j))$.
Finally, suppose that $\mathcal{G} \colon
\mathcal{M}_Y \to
\mathcal{M}_X\otimes\Lambda(t, dt)$ is a
right homotopy of DG
algebra maps from $\phi$ to $\psi$.  Then
$\QL(\mathcal{G})\circ G
\colon \L(W)_I \to \QL(\mathcal{M}_Y)$ is a
left homotopy of DG
Lie algebra maps from $\QL(\phi)\circ\rho_X$
to
$\QL(\psi)\circ\rho_X$.

We now come to the main point of the
appendix.

\begin{lemma}\label{lem:restricted DG Lie
homotopy}
Let $f \colon X \to Y$ be a map with a fixed
choice of Quillen
model $\mathcal{L}_f \colon \mathcal{L}_X
\to \mathcal{L}_Y$. Let
$A, B \colon S^n\times X \to Y$ be maps that
restrict to $A\circ
i_2 = B\circ i_2 = f \colon X \to Y$.
Suppose $\L(W)$ is the
Quillen minimal model of $X$, and $\L(W, v,
W')$ is the Quillen
model of $S^n\times X$ (see
\corref{cor:Quillen model S times X}).

\begin{enumerate}
\item[(i)] Suppose $A \sim B$ via a homotopy
$H$ relative to $X$, that is,
suppose the diagram (\ref{eq:unbased homotopy}) commutes.  Then
there exists a DG homotopy $\mathcal{H}
\colon \L(W, v, W')_I \to
\mathcal{L}_{Y}$ from $\mathcal{L}_A$ to
$\mathcal{L}_B$ that is
``relative to $\mathcal{L}_X$," in the sense
that the diagram
$$\xymatrix{\L(W)_I \ar[r]^-{p}
\ar[d]_{J_I}&
\L(W)
\ar[d]^{\mathcal{L}_f}\\
\L(W, v, W')_I \ar[r]_-{\mathcal{H}}&
\mathcal{L}_Y }
$$
(strictly) commutes.
\item[(ii)] Suppose further that $A$ and $B$
restrict to $A\circ
i_1 = B\circ i_1 = * \colon S^n \to Y$, and
the homotopy $H$ is
also relative to $S^n$, so that
(\ref{eq:based homotopy})
commutes.  Then there exists a DG homotopy
$\mathcal{H} \colon
\L(W, v, W')_I \to \mathcal{L}_{Y}$ from
$\mathcal{L}_A$ to
$\mathcal{L}_B$ that is ``relative to
$\mathcal{L}_{S^n\vee X}$,"
in the sense that the diagram (\ref{eq:based DG Lie homotopy})
(strictly) commutes.
\end{enumerate}
\end{lemma}

\begin{proof}
(i) In the Sullivan model setting, path
objects for
$\mathcal{M}_X$ and
$\mathcal{M}_{S^n}\otimes\mathcal{M}_X$ are
related as in the commutative diagrams
$$\xymatrix{\mathcal{M}_{S^n}\otimes\mathcal
{M}_X \ar[d]_{q_2} \ar[r]^-{j} &
\mathcal{M}_{S^n}\otimes\mathcal{M}_X\otimes
\Lambda(t,dt)\ar[r]^-
{p_i}\ar[d]^{q_2\otimes1}
& \mathcal{M}_{S^n}\otimes\mathcal{M}_X
\ar[d]_{q_2}\\
\mathcal{M}_X \ar[r]_-{j} &
\mathcal{M}_X\otimes\Lambda(t,dt)
\ar[r]_-{p_i} & \mathcal{M}_X }
$$
for $i = 0,1$. Applying $\QL$ to this
diagram, and following the
construction of the map $G$ as described
above, we obtain
homotopies $G_X \colon \L(W)_I \to
\QL(\mathcal{M}_X\otimes\Lambda(t,dt))$ and
$G_{S^n\times X}
\colon \L(W, v, W')_I \to
\QL(\mathcal{M}_{S^n}\otimes\mathcal{M}_X\otimes\Lambda(t,dt))$.
By constructing $G_{S^n\times X}$ so as to
extend $G_{X}$ on
$\L(W)_I$, we may assume these maps are
compatible, in the sense
that $\QL(q_2\otimes1)\circ G_X =
G_{S^n\times X}\circ J_I \colon
\L(W)_I \to
\QL(\mathcal{M}_{S^n}\otimes\mathcal{M}_X\otimes\Lambda(t,dt))$.
Combining this with the diagram obtained
from applying $\QL$ to
(\ref{eq:unbased Sullivan homotopy}), we
obtain a commutative
diagram
$$\xymatrix{\L(W)_I \ar[d]_{J_I} \ar[r]^-{G_X} &
\QL(\mathcal{M}_X\otimes\Lambda(t,dt))\ar[d]
^{\QL(q_2\otimes1)} \ar[r]^-{\QL(j)} &
\QL(\mathcal{M}_X)
\ar[d]^{\QL(\mathcal{M}_f)} \\
\L(W, v, W')_I \ar[r]_-{G_{S^n\times X}} &
\QL(\mathcal{M}_{S^n}\otimes\mathcal{M}_X\otimes\Lambda(t,dt))
\ar[r]_-{\QL(\mathcal{G})}
& \QL(\mathcal{M}_Y) }
$$
Furthermore, from the construction of $G_X$,
we have $\QL(j)\circ
G_{X}(sW, \widehat{W}) = 0$ so that
$\QL(j)\circ G_{X}=
\rho_X\circ p$.  It remains to check that we
can preserve the
properties of this diagram in lifting the
homotopies to
$\mathcal{L}_Y$.

So suppose that $\mathcal{L}_f$ is a given
Quillen model that
arises in the way described above from
$\phi_f \colon \L(W) \to
\mathcal{L}_Y\sqcup E(\QL(\mathcal{M}_Y))$.
With the above
ingredients, we obtain a diagram
$$\xymatrix@C=40pt{\L(W)_I \ar[d]_{J_I}
\ar[r]^-{\phi_f\circ p} &
\mathcal{L}_Y\sqcup E(\QL(\mathcal{M}_Y))
\ar@{->>}[d]^{\widehat{\rho_Y}}_{\simeq} \\
\L(W, v, W')_I \ar[r]_-
{\QL(\mathcal{G})\circ G_{S^n\times
X}}\ar@{.>}[ru]^{\mathcal{G}'} &
\QL(\mathcal{M}_{Y}), }
$$
which commutes since
$\widehat{\rho_Y}\circ\phi_f\circ p =
\QL(\mathcal{M}_{f})\circ\rho_X\circ p =
\QL(\mathcal{M}_{f})\circ\QL(j)\circ G_X$.
Therefore, we may lift
as indicated through the surjective quasi-isomorphism
$\widehat{\rho_Y}$ to obtain the homotopy
$\mathcal{G}'$ that
satisfies
$\widehat{\rho_Y}\circ\mathcal{G}' =
\QL(\mathcal{G})\circ G_{S^n\times X}$.
Finally, define
$\mathcal{H} = \pi\circ\mathcal{G}'$.  This
is a homotopy that
starts at a Quillen model for $A$, and ends
at one for $B$, since
we have
$$\begin{aligned}
\rho_Y\circ\mathcal{G}\circ\lambda_i &=
\rho_Y\circ\pi\circ\mathcal{G}'\circ\lambda_
i \\
&\sim
\widehat{\rho_Y}\circ\mathcal{G}'\circ\lambda_i =
\QL(\mathcal{G})\circ G_{S^n\times
X}\circ\lambda_i\\
 &= \QL(\mathcal{G})\circ
\QL(p_i)\circ\rho_{S^n\times X} =
\QL(p_i\circ\mathcal{G})\circ\rho_{S^n\times
X}
\end{aligned}$$
which equals
$\QL(\mathcal{A})\circ\rho_{S^n\times X}$
for $i = 0$
and $\QL(\mathcal{B})\circ\rho_{S^n\times
X}$ for $i = 1$. That
is, $\mathcal{G}\circ\lambda_0$ and
$\mathcal{G}\circ\lambda_1$
may be taken as Quillen models for $A$ and
$B$, respectively.

(ii) The argument for (i) is easily adapted
to establish (ii). In
this case, begin by constructing a homotopy
$$G_{S^n\vee X} \colon \L(W, v)_I \to
\QL(\mathcal{M}_X\otimes\Lambda(t,dt))\sqcup
\QL(\mathcal{M}_{S^n}\otimes\Lambda(t,dt))$$
so that $(j\sqcup j)\circ G_{S^n\vee X} =
(\rho_{X}\sqcup\rho_{S^n})\circ p \colon
\L(W, v)_I \to
\QL(\mathcal{M}_X)\sqcup\QL(\mathcal{M}_{S^n
})$.  Then extend
$\big( \QL(q_2\otimes1) \mid
\QL(q_1\otimes1)\big)\circ G_{S^n\vee
X}$ to $G_{S^n\times X} \colon \L(W, v,
W')_I \to
\QL(\mathcal{M}_{S^n}\otimes\mathcal{M}_X\otimes\Lambda(t,dt))$.
The argument now proceeds as before.
\end{proof}

Notice that in either case,
\lemref{lem:restricted DG Lie
homotopy} shows that we may choose a Quillen
model $\mathcal{L}_A$
of $A$ that restricts to $\mathcal{L}_f$ on
$\L(W)$ (the
restriction equals $\mathcal{L}_f$, and is
not just homotopic to
it).  This justifies a fact that we relied
upon for the definition
of the maps $\Phi$ and $\Psi$ in
\thmref{thm:isos Phi and Phi*}.
\end{appendix}

\providecommand{\bysame}{\leavevmode\hbox
to3em{\hrulefill}\thinspace}



\begin{thebibliography}{WL88b}

\bibitem[DZ79]{D-Z}
E.~Dror and A.~Zabrodsky, \emph{Unipotency
and nilpotency in
homotopy
  equivalences}, Topology \textbf{18}
(1979), no.~3, 187--197.

\bibitem[FH82]{F-H}
Y.~F{\'e}lix and S.~Halperin, \emph{Rational
{L}.-{S}. category
and its
  applications}, Trans.~Amer. Math.~Soc.
\textbf{273} (1982), no.~1, 1--38.

\bibitem[FHT01]{F-H-T}
Y.~F{\'e}lix, S.~Halperin, and J.-C. Thomas,
\emph{Rational
homotopy theory},
  Graduate Texts in Mathematics, vol. 205,
Springer-Verlag, New York, 2001.

\bibitem[Gan68]{Gan}
T.~Ganea, \emph{Cyclic homotopies}, Ill. J.
Math. \textbf{12}
(1968), 1--4.

\bibitem[Got65]{Got1}
D.~H. Gottlieb, \emph{A certain subgroup of the
fundamental group},
Amer. J. Math.
  \textbf{87} (1965), 840--856.

\bibitem[Got68]{Got2}
D.~H. Gottlieb, \emph{On fibre spaces and the
evaluation map}, Ann.
Math.
  \textbf{87} (1968), 42--55.

\bibitem[Got69]{Go1}
D.~H. Gottlieb, \emph{Evaluation subgroups
of homotopy groups},
Amer. J. Math.
  \textbf{91} (1969), 729--756.

\bibitem[HS79]{H-S}
S.~Halperin and J.~Stasheff,
\emph{Obstructions to homotopy
equivalences}, Adv.
  in Math. \textbf{32} (1979), 1233--279.

\bibitem[LS03]{L-S}
G.~Lupton and S.~B. Smith,
\emph{Rationalized evaluation subgroups
of a map: Sullivan models, derivations and
{$G$}-sequences}, Pre-print,
2003.

\bibitem[LW93]{L-W4}
K.~Y. Lee and M.~H. Woo, \emph{The {$G$}-sequence and the
{$\omega$}-homology
  of a {CW}-pair}, Topology Appl.
\textbf{52} (1993), no.~3, 221--236.

\bibitem[Opr95]{Oprea}
J.~Oprea, \emph{Gottlieb groups, group
actions, fixed points and
rational
  homotopy}, Lecture Notes Series, vol.~29,
Seoul National University Research
  Institute of Mathematics Global Analysis
Research Center, Seoul, 1995.

\bibitem[Spa89]{Spa89}
E.~H. Spanier, \emph{Algebraic {T}opology},
1st corrected
{S}pringer ed.,
  Springer-Verlag, New York, 1989.

\bibitem[SS]{Sch-St}
M.~Schlessinger and J.~Stasheff,
\emph{Deformation theory and
rational homotopy
  type}, preprint.

\bibitem[Sul78]{Su}
D.~Sullivan, \emph{Infinitesimal
computations in topology}, Inst.
Hautes Etudes
  Sci. Publ. Math. \textbf{47} (1978), 269--
331.

\bibitem[Tan83]{Tan83}
D.~Tanr{\'e}, \emph{Homotopie rationnelle:
mod\`eles de {C}hen,
{Q}uillen,
  {S}ullivan}, Lecture Notes in Mathematics,
vol. 1025, Springer-Verlag,
  Berlin, 1983.

\bibitem[WL88a]{L-W2}
M.~H. Woo and K.~Y. Lee, \emph{Homology and
generalized evaluation
subgroups of
  homotopy groups}, J. Korean Math. Soc.
\textbf{25} (1988), no.~2, 333--342.

\bibitem[WL88b]{L-W1}
M.~H. Woo and K.~Y. Lee, \emph{On the
relative evaluation
subgroups of a
  {CW}-pair}, J. Korean Math. Soc.
\textbf{25} (1988), no.~1, 149--160.

\end{thebibliography}

\end{document}